\documentclass[11pt]{amsart}
\usepackage{amssymb} 

\newtheorem{theorem}{Theorem}[section]
\newtheorem{lemma}{Lemma}[section]
\newtheorem{proposition}{Proposition}[section]

\theoremstyle{definition}
\newtheorem{definition}{Definition}[section]
\theoremstyle{remark}
\newtheorem{remark}{Remark}[section]
\numberwithin{equation}{section}
\newtheorem{corollary}{Corollary}[section]

\numberwithin{equation}{section}

\newcommand{\C}{\mathbb C}
\newcommand{\Z}{\mathbb Z}
\newcommand{\R}{\mathbb R}
\newcommand{\Q}{\mathbb Q}

\begin{document}

 \title[Critical points]{Global classification of isolated singularities
 in dimensions $(4,3)$ and $(8,5)$}
\author[L.Funar]{Louis Funar}
\address{Institut Fourier BP 74, UMR 5582, Universit\'e de
Grenoble, 38402 Saint-Martin-d'H\`eres cedex, France}
\email{funar@fourier.ujf-grenoble.fr}

\date{}

\begin{abstract} 
We characterize  those closed $2k$-manifolds admitting smooth maps 
into $(k+1)$-manifolds with only finitely many critical points, 
for $k\in\{2,4\}$. We compute then the minimal number of critical 
points of such smooth maps for $k=2$ and, under some fundamental group 
restrictions, also for $k=4$. 
The main ingredients are King's local classification of isolated 
singularities, decomposition theory, low dimensional cobordisms 
of spherical fibrations and 3-manifolds topology. 
\end{abstract}

\subjclass{57 R 45, 58 K 05, 57 R 60, 57 R 70}


\maketitle

\tableofcontents

\section{Introduction and statements}

Let $\varphi(M^m,N^n)$ denote the minimal number of critical points 
(not necessarily non-degenerate) 
of smooth  proper (i.e.  such that the inverse image of the 
boundary is the boundary) maps between the smooth manifolds $M^m$ and $N^n$. 
When superscripts are specified they denote the dimension of 
the respective manifolds. 

When we take 
$N^n=\R$ one obtains a classical homotopy invariant, namely 
the $F$-category of the manifold $M^m$. Using Morse-Smale theory Takens proved 
that the $F$-category  of a closed manifold is bounded by $m+1$ 
(see \cite{Takens2} for a refined upper bound in terms 
of the connectedness, the homology groups and the dimension). 
The $F$-category is bounded from below by the Lusternik-Schnirelman 
category. Related results for $h$-cobordisms were obtained by 
Pushkar and Rudyak (see \cite{PR}).

We are interested below in the case 
when $m\geq n\geq 2$ and the manifolds under consideration are compact 
and orientable. 
The invariant $\varphi$ is not anymore a homotopy invariant of 
the pair $(M^m,N^n)$ and actually it is sensitive to the smooth structures.  
There are actually only a few examples where this invariant can be 
calculated. The main problem in this area is to characterize those
pairs of manifolds for which $\varphi$ is finite non-zero
and then to compute its value (see \cite{CT1}, p.617). 
The problem of finding  non-trivial local isolated singularities 
was also stated by Milnor (\cite{Milnor}, section 11, p.100). 

In \cite{AF1} the authors found that,
in all codimension 
$0\leq m-n\leq 2$, if $\varphi(M^m,N^{n})$ is finite
then $\varphi(M^m,N^{n})\in\{0,1\}$, except for the
exceptional pairs of dimensions
$(m,n)\in\{(2,2), (4,3), (4,2)\}$. 
Further, if $m-n=3$ and there exists a smooth function  
$M^m\to N^n$ with finitely many critical points, all of them cone-like, 
then  $\varphi(M^m,N^{n})\in\{0,1\}$
except for the exceptional pairs of dimensions 
$(m,n)\in\{(5,2), (6,3), (8,5)\}$. 
Moreover, under the finiteness hypothesis,
$\varphi(M^m,N^n)=1$ if and only if $M^m$ is the connected sum 
of a  smooth fibration over $N^n$ with an exotic sphere and not a fibration 
itself. A short outline of this and related results for Lefschetz 
fibrations can be found in \cite{AF2}. 
These results comfort the idea that, in general, $\varphi$ is infinite. 
There is no enough room to plug inside a singular fiber enough of 
the topology of the manifold, at least when the codimension is much smaller  
with respect to the dimension. This might be compared to the Conjecture 
stated by Milnor in (\cite{Milnor}, p.100).

There are two essential ingredients in this result. First, 
the germs of smooth maps $\R^m\to \R^n$, in the given range of dimensions,  
having an isolated singularity at origin are locally topologically 
equivalent to a projection. Apparently the first who noticed that 
a smooth map of the $m$-disk in itself with a single interior critical point 
should actually be a homeomorphism for $m\geq 3$ 
was Heinz Hopf in the early thirties. The next step was taken by 
Antonelli, Church, Hemmingsen and Timourian who analyzed  in the 
seventies the local structure of smooth maps around a critical point, 
in codimension at most 2. Except a few exceptional dimensions
these maps are locally topologically equivalent to projections and 
thus they are topological fibrations. 
The second point is that  each critical point can be smoothed by 
removing a disk  and gluing it back, namely by adjoining a homotopy sphere. 
Thus all but possibly one critical points merge together.

The simplest exceptional case is that of pairs of  surfaces,
which admits an elementary treatment  and is
completely understood (see \cite{AF2} for explicit computations).
In \cite{FPZ} the authors described a series of examples  
in dimensions $(4,3)$, $(8,5)$ and $(16,9)$. In these dimensions 
one knew about the existence of smooth functions with genuine singular  
points since the papers of Antonelli (\cite{Ant1,Ant3}). 
These are particular Montgomery-Samelson fibrations 
(\cite{MonSam}) obtained from a spherical fibration by pinching a number 
of fibers to points. 
Around  such a critical point the map 
is locally topologically equivalent to the suspension over a 
Hopf fibration of spheres. This also explains the occurrences of dimensions 
$(4,3)$, $(8,5)$ and $(16,9)$. The novelty brought out in \cite{FPZ} is that 
the invariant $\varphi$ can be explicitly computed.  
There are obstructions to the clustering of  genuine singular points, 
whose number is determined by the algebraic topology of the manifold. 
In particular, one founds that $\varphi$ can take 
now arbitrarily large even values, contrasting with the situation 
treated before in \cite{AF1}.

The first aim of this paper is to explore the codimension 3 situation 
by 3-manifold topology methods. When the singularities are 
cone-like the general methods from \cite{AF1} permit to settle 
the problem, at least for $n\geq 4$. 
However, Takens (see \cite{Takens}) provided examples of  
wild (i.e. not cone-like) codimension 3 isolated singularities. 
In order to understand  wild singularities 
we revisit King's (\cite{King,King2}) local classification 
of isolated singularities in terms of 
topological data consisting of a local fibration and a 
vanishing compactum.  
We first obtain that there are actually uncountably many 
non-equivalent germs of  wild isolated singularities in codimension 3, 
so previous techniques cannot work.  The equivalence classes here 
correspond to right composition by (germs of) homeomorphisms.  
We introduce then the notion of removable singularity, 
where one enlarges the orbit of a singularity germ  
by allowing right compositions with limits of homeomorphisms.  
The key technical result is that codimension 3 singularities 
are removable, as soon as $n\geq 4$. In particular, we obtain: 

\begin{theorem}\label{local}
Let $M^{n+3}\to N^{n}$ be a smooth map with finitely many critical points. 
Then $M^{n+3}$ is diffeomorphic to $\Sigma^{n+3}\# Q^{n+3}$ where:  
\begin{enumerate}
\item $\Sigma^{n+3}$ is a homotopy sphere; 
\item  there is a smooth  (induced) map $Q^{n+3}\to N^n$ 
having only finitely many critical points and fulfilling:
\begin{enumerate}
\item if $n\geq 6$ or $n=4$ and the singularities were cone-like, 
then there are no critical points; 
\item  if $n=5$ then around  each critical point the map is 
smoothly equivalent to the cone over the  corresponding Hopf  fibration. 
\end{enumerate}
\end{enumerate}
\end{theorem}

We don't know whether wild singularities show up for general $(m,n)$.  
A more tractable problem in higher dimensions is to understand whether 
there exist smooth maps $M^m\to N^n$ with only cone-like isolated 
singularities. The local version of this problem consists 
of constructions of non-trivial germs of smooth (or even polynomial) 
maps and it  was actually reduced by Looijenga, 
Church and Lamotke (see \cite{Lo,CL}) to the existence of 
high dimensional fibered links. 
This was further pursued by Rudolph, Kauffman and Neumann (\cite{Rud,KN}) 
who considered general open book decompositions and more recently 
in \cite{AT}. 
The global question is how can we piece 
together conical extensions of open book decompositions 
to get well-defined maps  $M^m\to N^n$ and whether critical points could be 
cancelled.  

Milnor asked in 
(\cite{Milnor}, p.100) about the right  meaning 
of  being a {\em non-trivial} local 
isolated singularity and gave as an example the requirement 
that the local  generic fiber be a disk. It was shown in \cite{Lo} that an isolated singularity 
is trivial iff its germ is right equivalent to a topological projection 
(at least when the codimension is different from 4 and 5). 
In this paper we consider a weaker non-triviality condition, 
which we call removability,  whose meaning is that we can change the 
map in a small neighborhood of the critical point so that the 
new function has no (topological) 
critical points anymore. In other words the critical point does not contribute 
to $\varphi$. Notice that the non-trivial examples from 
\cite{CL} are actually removable isolated singularities. 
The search for non-removable global examples is therefore  rather subtle.

\begin{remark}
If the group of homeomorphisms of a contractible manifold $F$ 
(of dimension $m-n\geq 5$) which 
are identity on the boundary were contractible, then we 
could prove by the present methods that $\varphi(M^m,N^n)\in\{0,1\}$ whenever 
$n\geq m-n+1$, $n\geq 2$ and $(m,n)\not\in\{(4,3), (8,5), (16,9)\}$.   
However, the connectivity results from \cite{BLR} are not sufficient for 
obtaining the triviality of the fibrations with fiber $F$ when $n\geq m-n+1$. 
This is related to Milnor's Conjecture from 
(\cite{Milnor}, p.100) stating the triviality 
(here in the sense of removability) of local isolated 
cone-like singularities in this range of dimensions. 
This is similar to the constraints given in \cite{Rees} for the existence 
of even non-trivial cone-like isolated singularities.   
\end{remark}

The second aim of this paper is to use the local description of singularities 
in order to analyze the global topology of the involved manifolds 
in the case of exceptional dimensions. 
We  provide here a complete solution for this problem
in dimensions $(4,3)$ and $(8,5)$.  
We will show that the examples 
obtained in \cite{FPZ} are essentially 
all examples with finite $\varphi$, up to a certain operation  
on manifold fibrations called {\em fiber sum}.

The case of dimensions $(4,3)$ still relies on the local classification 
of singularities due to Church and Timourian but 
the case of dimensions $(8,5)$ needs 
Theorem \ref{local}, which enables us to work with local models which are 
cones over Hopf fibrations. 
The main difficulty is then to determine which are the smooth manifolds 
obtained by piecing together these local models. 

It appears that the gluing process is surprisingly rigid 
and we obtain always connected sums of $S^2\times S^2$ and respectively 
of $S^4\times S^4$. 
Eventually, we can also glue along the fibers a multi-spherical 
fibration which does not contribute to the critical points. 
The result of such a construction is abstracted below:

\begin{definition}\label{genfsum0}
Assume given the following data: 
\begin{enumerate}
\item  an $F$-bundle $g:E^m\to Q^{k+1}$, and a smooth 
map $h:X^m\to D^{k+1}$ with 
finitely many critical points  such that the generic fiber of $h$ is  
$F$ and 
\item an isomorphism between the 
$F$-fibrations $g|_{\partial E}:\partial E\to \partial Q$ 
and $h|_{\partial X}:\partial X \to \partial D^{k+1}$, 
(thus  $\partial Q^{k+1}=S^k$).
\end{enumerate} 
Glue together $g$ and $h$  to obtain a smooth map 
$f:E^m \cup X^m \to Q^{k+1}\cup D^{k+1}$, with finitely many critical points, 
by identifying the boundary fibrations by means of the fixed  bundle 
isomorphism. 
We call $f$ the  {\em generalized fiber sum} of $g$ and $h$ and denote it by 
$g\oplus h$; sometimes we mention also the isomorphism type $\alpha$  
of the boundary $F$-fibration  $h|_{\partial X}$ over $S^k$ and write 
$g\oplus_{\alpha} h$. Notice that this fibration might well be non-trivial, 
in contrast to the usual fiber sum construction used in \cite{FPZ}. 
One says that $h$ is the {\em non-trivial summand} and 
$g$ the {\em fibrewise summand} of the fiber sum.

When $g:E^m\to Q^{k+1}$ and $h:X^m\to D^{k+1}$ are obtained from  
the smooth maps $G:\widehat{E}^m\to \widehat{Q}^{k+1}$,
$H:\widehat{X}^m\to S^{k+1}$ between closed manifolds,  
by removing trivial fibrations over disks, then  $g\oplus h$ is 
the {\em fiber sum} of the maps $G$ and $H$. 
\end{definition}

The main result of this paper is the following structure Theorem:

\begin{theorem}\label{complete}
Let $M^{2k}$ and $N^{k+1}$ be closed orientable manifolds 
having finite $\varphi(M^{2k},N^{k+1})$, where   $k\in\{2,4\}$. Then 
\begin{enumerate}
\item either $M^{2k}$ is diffeomorphic to $W^{2k}\#\Sigma^{2k}$, 
where 
$W^{2k}$ is a fibration over $N^{k+1}$ and 
$\Sigma^{2k}$ is a homotopy $2k$-sphere. 
In this case 
$\varphi(M^{2k},N^{k+1})\in\{0,1\}$. 
\item or else, $M^{2k}$ is diffeomorphic to the iterated fiber sum 
$W^{2k}\oplus_{j=1}^D\#_{r_j} S^{k}\times S^{k} 
\# \Sigma^{2k}$, with $r_i\geq -1$, but at least one $r_i\geq 0$, 
where: 
\begin{enumerate}
\item the fibrewise summand $W^{2k}$ is a $S^{k-1}$-fibration over 
some closed manifold $\widehat{N^{k+1}}$, where 
$\widehat{N^{k+1}}\to N^{k+1}$ is a non-ramified 
covering of  finite degree $D$.  
\item the fiber sum $W^{2k}\oplus_{j=1}^D\#_{r_j} S^{k}\times S^{k}$  is the fiber sum along the fiber 
of  the fibration $W^{2k}\to N^{k+1}$ (which 
consists of $D$  disjoint copies of $S^{k-1}$) 
with the disjoint union   
$\sqcup_{j=1}^D\#_{r_j} S^{k}\times S^{k}$. 
Each connected component of the later 
is endowed with a smooth map 
$\#_{r_j} S^{k}\times S^{k}\to S^{k+1}$
whose generic fiber is $S^{k-1}$ (to be defined below).  
\item   $\#_{-1} S^{k}\times S^{k}$ denotes 
$S^{k-1}\times S^{k+1}$, 
$\#_{0} S^{k}\times S^{k}$ denotes $S^{2k}$ and  
$\Sigma^{2k}$ a homotopy $2k$-sphere.  
\end{enumerate}
Moreover, in this case we have:  
\[\varphi(M^{2k}, N^{k+1})\leq  2r_1+\cdots +2r_D+2D\] 
\end{enumerate}
\end{theorem}

\begin{remark}
The fiber sum  of $W^{2k}$ 
with $\#_{-1} S^{k}\times S^{k}=S^{k-1}\times S^{k+1}$ has no effect 
on the topology of $W^{2k}$ since it 
is diffeomorphic to $W^{2k}$. Thus, if all $r_i=-1$ then we recover the 
manifold $W^{2k}$ and thus $M^{2k}$ fibers over $N^{k+1}$. 
\end{remark}

\begin{remark}
There is one  extra piece in the fiber sum data, namely the choice 
of a bundle isomorphism between  each pair of trivial 
$S^{k-1}$ fibrations over $S^k$. Such bundle isomorphisms up to isotopy 
correspond to the set of homotopy classes of maps 
$S^k\to {\rm Homeo}^+(S^{k-1})$. 

If $k=2$ this set is in bijection with $\pi_2(SO(2))=0$ and thus the 
fiber sum construction is unambiguously defined. 

When $k=4$, using Hatcher's solution 
to the Smale Conjecture, this set is in bijection with  
$\pi_4({\rm Homeo}^+(S^{3}))=\pi_4(SO(4))=\Z/2\Z+\Z/2\Z$ and 
the fiber sum depends on the choice of $D$ elements of this 
group, which will be called gluing parameters in the sequel. 
\end{remark}

When $M^{2k}$ satisfies the first alternative of Theorem \ref{complete} 
it is a topological fibration over $N^{k+1}$. If $M^{2k}$ is not a 
smooth fibration then $\varphi(M^{2k},N^{k+1})=1$. 

\begin{remark}
The manifold $M^m=\Sigma^m\# S^{m-n-1}\times S^{n+1}$ is not diffeomorphic 
to $S^{m-n-1}\times S^{n+1}$ if  $\Sigma^m$ is an exotic sphere (see  
\cite{Sch}). This yields effective examples 
where $\varphi=1$.  For instance, if $\Sigma^8$ is the exotic 8-sphere 
which generates the group $\Gamma_8=\Z/2\Z$ of homotopy spheres of dimension 
8 then    
$\varphi(\Sigma^8\# S^3\times S^5,S^5)=1$ (see e.g. \cite{FPZ}).
\end{remark}

When the second alternative of Theorem \ref{complete} holds,  
$M^{2k}$ will be called {\em non-fibered} (over $N^{k+1}$). 
Of course these furnish the most interesting examples with non-trivial 
$\varphi$. Moreover, Theorem \ref{complete} leads to 
effective topological obstructions to the finiteness of $\varphi$, as follows:

\begin{corollary}\label{gp}
\begin{enumerate}
\item If $\varphi(M^8,N^5)$ is finite and  $M^8$ is  non-fibered over $N^5$  
then $\pi_1(M^8)$ is a finite index subgroup of $\pi_1(N^5)$.
\item If  $\varphi(M^4,N^3)$ is finite and $M^4$ is  non-fibered over $N^3$  
then $\pi_1(M^4)\cong \pi_1(N^3)$. 
\end{enumerate}
\end{corollary}

Further,  there are global constraints 
of topological nature to the clustering of 
genuine critical points. This situation seems rather exceptional and 
specific to the dimensions $(4,3)$, $(8,5)$ and $(16,9)$. 
Moreover, the number of genuine critical points is independent on the 
smooth function considered and coincides with $\varphi$, which is 
therefore controlled by the topology, as follows:

\begin{theorem}\label{compute}
Let $M^{2k}$ and $N^{k+1}$  
be closed orientable manifolds with finite 
$\varphi(M^{2k},N^{k+1})$, with $k\in\{2,4\}$, 
and $M^{2k}$ non-fibered over $N^{k+1}$. Let  $r=r_1+r_2+\cdots+r_D$, 
where $r_j, D$ are those furnished by Theorem \ref{complete}.
\begin{enumerate} 
\item If $k=4$, assume that  $\pi_1(M^8)\cong \pi_1(N^5)$ is a co-Hopfian group 
or a finitely generated free non-abelian group. Then 
$\varphi(M^8,N^5)=2r+2D$.    
\item If $k=2$,  then $\varphi(M^4,N^3)=2r+2D$. 
\end{enumerate}
\end{theorem}

One might conjecture that the assumptions on $\pi_1(M^8)$,  
when  $k=4$, are superfluous and that 
we have $\varphi(M^8,N^5)=2r+2D$  
for any non-fibered $M^8$, so that the right hand side 
is independent on the particular choices we made for 
writing the non-fibered $M^8$ as a fiber sum.

\begin{remark}
If the fibration of $W^{2k}$ is the connected sum of  $S^1$-products 
of Hopf fibrations, namely  it is the fibration $\#_c S^1\times S^{2k-1}\to \#_c S^1\times S^k$ 
then 
the manifold $M^{2k}$ is diffeomorphic to 
$\Sigma^{2k}\#_{r+c} S^k\times S^k \#_c S^1\times S^{2k-1}$ and 
the value for $\varphi$ given above is consistent with that obtained in  
the main Theorem of \cite{FPZ}. 
\end{remark}

\begin{remark}
For $k\in \{2,4\}$, the  isomorphism classes of smooth 
$S^{k-1}$-fibrations over $\#_c S^1\times S^k$ are classified by 
elements of $\pi_{k-1}(SO(k))^c$ which are either 
$(\Z\oplus\Z)^c$ (for $k=4$) or else $\Z^c$ (for $k=2$). 
\end{remark}

As a consequence of the present paper and \cite{AF1} we obtain the 
global classification of isolated singularities in codimension 
at most 3 except for the dimensions 
$(4,2)$, $(5,2)$ and $(6,3)$. 

The local classification of isolated cone-like singularities 
in dimensions $(6,3)$ corresponds to that of fibered links 
$\sqcup_{j=1}^k S^2$ into $S^5$ for $k\geq 2$, which 
is given by the linking matrix. However, new methods are necessary 
for obtaining the global classification  of isolated 
singularities in dimensions $(6,3)$. For instance,  
we don't know whether in the process of pasting together 
local singularities we could have any cancellation of critical points.  

The case of dimensions $(4,2)$ offers another interesting perspective, 
because of the abundance of non-equivalent fibered links in $S^3$. 
Their suspensions yield examples in dimensions $(5,2)$.

The plan of the paper is as follows. 
The first section describes the local classification of 
isolated singularities of codimension at most 2, following Church and 
Timourian and of cone-like isolated singularities of codimension 3 and 
dimension at least 4. The next section deals with wild isolated singularities. 
In order to understand them better we revisit King's classification of 
cone-like and general singularities. Each local singularity corresponds to 
some topological data consisting of a fibration $E$ 
(like in a open book structure) over a sphere, a vanishing compactum 
${\mathcal A}\subset E$ whose complementary is a trivial fibration 
over the sphere and a retraction $r$ of $E$ onto the singular fiber $V$ 
obtained from the regular fiber by crushing the vanishing compactum. 
King's classification of isolated singularities states that this 
data up to a certain equivalence relation determines the right equivalence 
class of the germ. In the cone-like case  the statement is more precise 
as one showed that all such data with tame vanishing compactum ${\mathcal A}$ 
are realized by a cone-like singularity.  
Our contribution is to explain that in general, what we require  is that  
the mapping cylinder of the retraction $r$ be a topological manifold. 
We intend further to use this to the classification of 
wild isolated singularities of codimension 3. It appears that 
fibrations whose fiber is a manifold with boundary of dimension at most 3  
and whose boundary is already a product (over a high dimensional sphere) 
tend to be trivial. This is a consequence of the contractibility 
of the group of homeomorphisms, due to Earle, Eells and Schatz for 
large surfaces and to Hatcher for irreducible 3-manifolds. Some work is needed 
to understand that the reducible case would lead to non-trivial 
homotopy of the singular fiber, which cannot happen in large dimensions.  
When the associated fibration is trivial the mapping cylinder of the retraction 
$r$ is the quotient of a product by the vanishing compactum. 
Decomposition theory will tell us that such a quotient space is a manifold 
only when the crushed compactum is cellular and then the quotient surjection 
is a near-homeomorphism. This actually explains how the 
wild singularities can be constructed (by using suitable Artin -Fox wild arcs 
in $S^3$)  and why their local structure is almost-trivial: although 
there are uncountably many wild singularities  non-equivalent 
by right composition with homeomorphisms, they are all equivalent 
by right composition with a near-homeomorphism. 
The near-homeomorphisms are precisely the cell-like maps, as 
it was proved by Siebenmann and thus they are limits of homeomorphism.
This leads to Theorem \ref{local}.  
We further want to describe a global structure result for 
manifolds $M^{2k}$ admitting smooth maps with finitely many 
critical points into $N^{k+1}$, for $k\in\{2,4\}$, where we know 
how the local germs look like. We split off a fibration factor to  obtain 
a sub-manifold $X^{2k}$ having  finite $\varphi(X^{2k},D^{k+1})$, which is  
called a  {\em disk block}. If $\varphi(X^{2k},S^{k+1})$ is finite then 
$X^{2k}$ is called a {\em spherical block}. 
We can first use general 
classification results of manifolds to  
find the diffeomorphism type of spherical blocks in dimension 8 and of 
the homeomorphism type of spherical blocks in dimension 4.  
These manifolds are  actually Montgomery-Samelson 
fibrations with finite branch set and their algebraic invariants are 
known since the seventies. In order to pursue further we need 
some ad-hoc techniques. 
We remark that disk blocks are obtained by gluing together the 
cones over the Hopf fibrations. Constructing such a disk block amounts to find 
a cobordism between the associated fibrations. 
However the theory of cobordisms of spherical fibrations 
is quite simple in these dimensions. This shows that such disk blocks are 
connected sums of spherical blocks from which one delete a neighborhood of 
a generic fiber. This gives the structure Theorem \ref{complete} in dimension 
8. In order to obtain the 4-dimensional situation we have to understand 
the spherical blocks with 2 critical points. As a consequence of Cerf's Theorem 
$\Gamma_4=0$ these are smooth $S^4$. We can further show that all 
spherical blocks are connected sums of $S^2\times S^2$ up to adding some 
homotopy 4-sphere. 
We obtain also an upper bound for the number $\varphi$ 
of critical points, which we are able in the last section  to prove that it 
is sharp. In fact, when the fundamental group is co-Hopfian 
we can show that $\varphi$ is determined by the homology of the 
manifolds involved. 
Using Thurston's geometrization Conjecture 
in dimension 3 (as settled by Perelman) 
and the results of Wang, Yu and Wu which characterize 
all 3-manifolds whose fundamental group is non co-Hopfian we can settle 
completely the case of dimensions $(4,3)$, by showing that the previous 
value of $\varphi$ is independent on the choices of intermediary coverings.

\vspace{0.2cm}

{\em Proviso}.  We will make use through out the present paper 
of the validity of the Thurston geometrization Conjecture and, in particular, 
of the Poincar\'e Conjecture in dimension 3, as settled by Perelman 
(see \cite{BBB} for a detailed proof), 
but we will specify its use each time when applied. 

\vspace{0.2cm}

{\bf Acknowledgements}. The author is indebted to Christine Lescop, 
Patrick Popescu-Pampu, Colin Rourke, Mihai Tib\u{a}r and Ping Zhang for useful discussions on 
this topic and acknowledge partial support from the 
ANR Repsurf: ANR-06-BLAN-0311.

\section{Local structure of isolated singularities 
of small codimension}


\subsection{Local classification of isolated singularities of codimension at most 2}

Let $f:M^m\to N^{n}$ be a smooth map with finitely many critical points 
throughout this section.

Church and Timourian (see \cite{CT1,CT2}) proved the following: 

\begin{proposition}\label{dim}
If  $m-n\leq 2$ and $x\in M$ then $f$ is locally topologically
equivalent to one of the following local models:
\begin{enumerate}
\item a projection;
\item the map $g:\C\to \C$ given by $g(z)=z^d$, $d\in \Z_+$, when $m=n=2$;
\item the Kuiper map $\tau_{2,1}:\C\times \C \to \C \times \R$ given by
$\tau_{2,1}(z,w)=(2z\bar{w}, |w|^2-|z|^2)$, when $m=4, n=3$;
\item a map $\rho:\R^4\to \R^2$ which is locally topologically equivalent
to the projection except at one point, when $m=4$ and $n=2$.
\end{enumerate}
\end{proposition}
\begin{remark}\label{kuiper}
There are more general Kuiper maps (see \cite{Milnor}, p.103)
$\tau_{a,k}:\R^{2ak}\to \R^{a+1}$, where $a\in \{2,4,8\}$ and $k\in \Z_+$, defined  as 
\[ \tau_{a,k}: A^k\times A^k\to A\times \R, \, 
\tau_{a,k}(x,y)=(\langle x, y\rangle, |x|^2-|y|^2)\]
Here $A\in\{\C, \mathbb H, \mathbb O\}$, where $\mathbb H$ denotes the 
quaternions, $\mathbb O$ denotes the Cayley numbers, 
$a\in\{2,4,8\}$   corresponds to the dimension of $A$ as a real vector 
space and  $\langle *, *\rangle$ denotes a Hermitian inner product. 
The map $\tau_{a,1}$ restricts to the Hopf fibration between the  corresponding unit spheres 
$S^{2a-1}\to S^{a}$ and  thus  it is equivalent to the 
suspension of the Hopf fibration.  For any $k\in\Z_+$ the map 
$\tau_{a,k}$ has an isolated singularity at the origin whose local 
fiber is a disk bundle over $S^{ak-1}$. 
\end{remark}

One might wonder whether there are only finitely many non-equivalent 
germs of isolated singularities, when $m-n=3$.  
First we should add one more germ  appearing  
when $m=8$ and $k=5$, namely the Kuiper map $\tau_{4,1}$, 
which is the cone  of the Hopf fibration $S^7\to S^4$.
However, it is still not enough (see the next section) 
since there exist infinitely many non-equivalent 
isolated germs, in general. But there are some slightly weaker 
statements which are true. We could restrict ourselves to 
cone-like isolated singularities, as those arising from functions which are 
locally real analytic around the critical points.  
On the other hand we will define the concept of removability and show that 
those singularities which are not cone-like are actually removable.

The aim of the next two sections is to give a more 
conceptual and self-contained proof of Proposition \ref{dim} 
 and its various extensions to codimension 3. 

\begin{definition}
Let $V=f^{-1}(f(x))$, where $x$ is a critical point. 
Following King (see \cite{King})  the singular point $x$ is called {\em cone-like} if it admits a  
cone neighborhood in $V$, i.e. 
there exists some closed manifold $L\subset V\setminus\{x\}$ and a neighborhood 
$N$ of $x$ in $V$ which is homeomorphic to  the cone $C(L)$ over $L$. Recall that the  
cone is defined as  the quotient $C(L)=L\times (0,1]/L\times \{1\}$. 
Then the manifold $L$ is called the local link at $x$.  If $x$ is not cone-like then $x$ (and  also 
$V$) are  called {\em wild}. An isolated point will be 
considered to be cone-like, as being the cone over the empty set. 
\end{definition}

The first examples of smooth maps with isolated wild singularities were obtained 
by Takens (see \cite{Takens}) in codimension 3.  
From \cite{CT1,CT2} such examples  cannot occur in smaller codimension 
because of the following:  

\begin{proposition}\label{codim2}
Isolated singularities of smooth functions in codimension  at most 2 are cone-like. 
\end{proposition}
\begin{proof}
This is clear from Proposition \ref{dim}, 
except when $m=4$ and $k=2$.  This case is settled 
in (\cite{CT2}, Lemmas 2.1 and 2.4).  
We will give another proof in the next section. 
\end{proof}

\subsection{Local classification of cone-like isolated singularities in codimension 3}
 
Denote by $S(f)$ the {\em topological branch locus} of  the map $f$. Namely 
$x\not\in S(f)$ if $f$ is locally topologically equivalent to a projection 
around $x$.

\begin{proposition} 
Let $f:M^{n+3}\to D^n$, $n\geq 4$ be a smooth proper map from the compact 
manifold $M^{n+3}$ onto the $n$-disk  having a single critical point 
$x\in M$ and suppose that $f(x)=0$. Then one of the following holds:
\begin{enumerate}
\item  $x\not\in S(f)$; 
\item  $x$ is an isolated point of $V=f^{-1}(0)$; 
\item or else $x\in S(f)$, $V$ is wild at $x$. 
\end{enumerate}
\end{proposition}

In the proof of this Proposition we use the 
following lemma stated also in \cite{FPZ}: 
\begin{lemma}\label{connect}
 If $M^{n+q}$ and $N^{n}$ are smooth manifolds and
$f:M^{n+q}\rightarrow N^n$ is a smooth map with finitely many critical
points, then the inclusions $M\setminus
V\hookrightarrow M$ and $N\setminus B\hookrightarrow N$ are $(n-1)$-connected.
\end{lemma}
\begin{proof}
For the sake of completeness we give here the proof. 
The result is obvious for  $N\setminus
B\hookrightarrow N$. It remains to prove
that $\pi_j(M,M\setminus V)\cong
0$, for $j\leq n-1$. Take $\alpha:(D^{j+1},
S^{j})\rightarrow (M, M\setminus V)$  to be an arbitrary
smooth map of pairs.
Since the critical set $C(f)$ of $f$ is finite and included in
$V$, there exists a small homotopy of $\alpha$ relative to
the boundary such that the image $\alpha(D^{j+1})$ avoids $C(f)$.
By compactness there exists a neighborhood $U$ of $C(f)$ 
consisting of disjoint balls centered at the critical points such that 
$\alpha(D^{j+1})\subset M\setminus U$.  We can arrange by a small isotopy 
that $V$ become transversal to $\partial U$.

Observe further that $V\setminus U$ consists of regular
points of $f$ and thus it is a properly embedded sub-manifold of
$M\setminus U$. General transversality arguments show that $\alpha$
can be made transverse to  $V\setminus U$ by a small
homotopy. By dimension counting this means
that $\alpha(D^{j+1})\subset M\setminus U$ is disjoint from $V$ and thus
the class of $\alpha$ in  $\pi_j(M,M\setminus V)$ vanishes.
\end{proof}

\begin{lemma}\label{fiber}
If $M^{n+3}$ is contractible  and $n\geq 4$ then the 
fiber $F^3$ of  the fibration $f:M^{n+3}\setminus V\to D^n\setminus\{0\}$ is 
diffeomorphic to a 3-disk. 
\end{lemma}
\begin{proof}
Lemma \ref{connect} states that $\pi_j(M,M\setminus V)=0$, for $j\leq n-1$ and 
thus $\pi_j(M\setminus V)=0$, for $j\leq n-2$. 
The long exact sequence in homotopy of the 
fibration $f:M^{n+3}\setminus V\to D^n\setminus\{0\}$ implies that 
the  generic fiber $F^3$ has $\pi_j(F)=0$, for $j\leq n-3$. 
If $n\geq 5$ then $F^3$ is  contractible and thus, by the Poincar\'e Conjecture, 
diffeomorphic to a 3-disk.

If $n=4$ then $F^3$ is a simply connected 3-manifold with boundary and hence 
its boundary consists of several 2-spheres. Thus $H_2(F)$ is free abelian. 
The base space of the fibration $f:M\setminus V\to D^4\setminus\{0\}$ 
is homotopy equivalent to $S^{3}$ and thus there is a Wang exact sequence 
for this fibration which reads:  
\[ 0=H_3(F)\to H_3(M\setminus V)\to H_0(F)\to H_2(F)\to H_2(M\setminus V)=0\]
Assume that $\pi_3(M\setminus V)\neq 0$. Then, using Hurewicz 
the group $H_3(M\setminus V)\neq 0$ is a 
non-trivial subgroup of $H_0(F)=\Z$ and hence it is isomorphic to $\Z$.  
This implies that the rank of $H_2(F)$  over $\mathbb Q$ is zero and hence 
$H_2(F)=0$ because it is free abelian. By Hurewicz, we have 
$\pi_2(F)=0$ and thus, by the Poincar\'e Conjecture, 
$F$ is diffeomorphic to the 3-disk, as claimed. 

Otherwise, $\pi_3(M\setminus V)=H_3(M\setminus V)=0$ 
and the Wang sequence above gives us $H_2(F)=\Z$ and 
$H_4(M\setminus V)=0$.  By the Poincar\'e Conjecture a simple connected 
3-manifold is a holed disk and thus $F=S^2\times [0,1]$. 
Moreover, from Hurewicz we obtain 
$\pi_4(M\setminus V)=0$. Furthermore, from the long exact sequence in homotopy 
of the fibration  $f:M\setminus V\to D^4\setminus\{0\}$ we derive: 
\[ 0=\pi_4(M\setminus V)\to \pi_4(D^4\setminus\{0\})\to \pi_3(F)=\pi_3(S^2)\to \pi_3(M\setminus V)=0\]
But $\pi_4(S^3)=\Z/2\Z$ while $\pi_3(S^2)=\Z$, which contradicts the exacteness 
of the sequence above. This proves the claim. 
\end{proof}

\begin{proposition}
If the isolated critical point $x$ of $f:M^{n+q}\to N^n$ is cone-like then there exists 
a compact manifold with boundary $K^{n+q}\subset M^{n+q}$ containing $x$ 
and a disk neighborhood $D^n$ containing the critical value $f(x)$ such that: 
\begin{enumerate} 
\item the restriction $f:K^{n+q}\to D^n$ is proper and has only one critical point $x$; 
\item $K^{n+q}$ is contractible. 
\end{enumerate}
\end{proposition}
\begin{proof}
This is basically contained in \cite{King}. Here is a simple proof. 
We consider a small enough manifold neighborhood $K^{n+q}$ of $x$ intersecting transversely 
$f^{-1}(D^n)$ and containing only one critical point. 
Let $F^q$ be the fiber. Then $K^{n+q}$  deformation retracts onto 
$V\cap K^{n+q}$ which 
is a cone $C(L)$ and thus contractible. 
\end{proof}

\begin{remark}
Actually if $n+q\neq 4,5$ we can assume that $K^{n+q}$ is homeomorphic to a ball and 
its boundary to a sphere (see again \cite{King}). 
This  implies that the fiber $F^3$ is simply connected also when $n=3$.   
\end{remark}

\begin{proposition}\label{dim3}
If  $m-n=3$, $n\geq 4$ and $x\in M^{m}$  is a cone-like isolated 
singularity then $f$ is locally topologically
equivalent to one of the following local models:
\begin{enumerate}
\item a projection;
\item the Kuiper map $\tau_{4,1}:\mathbb H\times \mathbb H \to \mathbb H \times \R$, 
when $m=8, n=5$ (see Remark \ref{kuiper});
\end{enumerate}
\end{proposition}
\begin{proof}
If $x$ is an isolated point of $V$ then, by 
\cite{Tim}, $f$ is locally topologically equivalent to the cone of a 
Hopf fibration and thus to the Kuiper map $\tau_{4,1}$. 
Henceforth, we assume that $z$ is not 
an isolated point in $f^{-1}(f(z))$, for any $z\in M$. 
Then $f$ has a cone-like critical point $x$ implies that $x\not\in S(f)$ 
because the local fiber is a disk and 
hence the result follows from \cite{CL}. 
\end{proof}


\subsection{King's classification revisited}
We consider germs of isolated singularities up to right composition by a homeomorphism,  
as described  by King  in  \cite{King} for continuous functions 
 and in \cite{King2} for smooth functions. 

Let $f:N_1^{n+q}\to N_2^n$ be a smooth map with finitely many 
critical points. Fix a critical point $p$. 
A smooth codimension zero compact  
manifold  with boundary $M^{n+q}\subset N_1^{n+q}$ containing $x$ in 
its interior such that the restriction $f:M^{n+q}\to D^n$ is a proper smooth 
map having only one critical point 
onto the disk $D^n\subset N_2^n$ is called an {\em adapted} neighborhood of the 
critical point. According to the construction of \cite{King} 
(see below) there exist arbitrarily small adapted neighborhoods of $p$. 

Let $f^{-1}(S^{n-1})=E^{n+q-1}\subset \partial M^{n+q}$, $V=f^{-1}(f(p))\subset M^{n+q}$ and 
$F^q$ denote the generic fiber, which is a manifold with non-empty boundary.  

The restriction $f:E^{n+q-1}\to S^{n-1}$ is a fibration with fiber $F^{q}$. 
Moreover, the restriction to the boundary is a trivial fibration 
$\partial E \to S^{n-1}$.  
Further $\partial V\cong \partial F$ is a closed submanifold in $\partial M$ 
having trivial normal bundle and which is endowed with a natural trivialization 
$t:\partial F\times D^n\to \partial M$, extending the trivial fibration 
$\partial E\to S^{n-1}$.  
Specifically we obtain $t$ by using the identification of 
$f:\partial M \setminus {\rm int}(E) \to D^n$ with the  trivial  
fibration over $D^n$. 

There is a retraction $r:E\to V$ obtained by contracting a 
vanishing compact ${\mathcal A}\subset E$ as follows: 

\begin{enumerate}
\item  $E\setminus {\mathcal A}$ is a product, i.e. it is homeomorphic to  
$W \times S^{n-1}$ such that the restriction of $f:E\to S^{n-1}$ to 
$E\setminus {\mathcal A}$, namely 
$f|_{E\setminus {\mathcal A}} : E\setminus {\mathcal A} \cong  W \times S^{n-1}\to S^{n-1}$ is the projection 
on the second factor.  In particular, the boundary 
fibration $f|_{\partial E} : \partial E\to S^{n-1}$ is trivial; 
\item $W=F\setminus A$, where $A\subset F$ is some compact and ${\mathcal A}\cap F =A$.  
Here $F$ is an arbitrary fiber in $E$ and a priori $A$ might depend on the particular fiber; 
\item $V= F/A$ and the class of $A$  in $F/A$ is the singular point $p$ of the fiber $V$. 
Thus $V\setminus\{p\}$ is homeomorphic to $W=F\setminus A$;  
\item the retraction $r:E\to V$ is obtained as follows: 
$r|_{E\setminus {\mathcal A} }: E\setminus {\mathcal A} \cong  W \times S^{n-1}\to W\cong V-\{p\}$
is given by the projection on the first factor of the product 
and $r(x)=p$, if $x\in {\mathcal A}$. 
\end{enumerate}

Let $\mathfrak{E}_{n,q}$ be the set of data 
$(E\to S^{n-1},  {\mathcal A}\subset E)$ as above up to the following  
equivalence relation. We set 
$(E_1\to  S^{n-1},{\mathcal A}_1)\sim (E_2\to S^{n-1},{\mathcal A}_2)$ 
if there exist sub-fibrations $E_j'\subset E_j$ 
such that ${\mathcal A}_j\subset E_j'$, 
which are homeomorphic  by a fibered  
homeomorphism $\varphi:E_1'\to E_2'$ (over $S^{n-1}$)
and there exists a fibered 
isotopy $h_t:E_1'\setminus {\mathcal A}_1\to E_2'$ such that 
\begin{enumerate}
\item $h_1$ is the restriction of $\varphi|_{E_1'\setminus {\mathcal A}_1}$; 
\item $h_0(E_1'\setminus {\mathcal A}_1)=E_2'\setminus {\mathcal A}_2$ 
and the fibered homeomorphism 
$h_0:E_1'\setminus {\mathcal A}_1\to E_2'\setminus {\mathcal A}_2$ 
comes from a homeomorphism between the compactified fibers i.e. 
the fiber restriction 
$h_0:F_1-A_1 \to F_2-A_2$ is induced from a 
homeomorphism $F_1/A_1 \to F_2/A_2$. 
\end{enumerate}

\begin{proposition}[King]
The class  $\delta(f)=(E\to S^{n-1},  {\mathcal A}\subset E)$ in 
$\mathfrak{E}_{n,q}$  is a complete invariant of  the germ of $f$, up to right 
equivalence i.e. up to right composition by a germ of  homeomorphism of $M$ at $p$. 
\end{proposition}
\begin{proof}
This is the main result of \cite{King} in a slightly changed setting. 
\end{proof}

\begin{definition}
The compact $A\subset F$ is {\em tame} if it has a mapping cylinder neighborhood i.e. 
a neighborhood $N$ which is a codimension zero compact submanifold of $F$ 
which deformation retracts onto $A$ such that $N\setminus A$ is homeomorphic to 
$\partial N \times (0,1]$.  
\end{definition}

\begin{remark}
If $A$ is a polyhedron in $F$ then it has a regular neighborhood and thus it is tame. 
If $A$ is an Artin-Fox wild arc (see \cite{ArF}) in $D^3$ then it is not tame. 
\end{remark}

\begin{proposition}
Let $\mathfrak{E}_{n,q}^{PL}\subset \mathfrak{E}_{n,q}$ be the subset of equivalence 
classes of pairs  $(E\to S^{n-1},  {\mathcal A}\subset E)$ 
for which ${\mathcal A}$ is tame. Then $\delta(f)\in \mathfrak{E}_{n,q}^{PL}$ classifies 
isolated  cone-like singularities of germs up to right equivalence.  
\end{proposition} 
\begin{proof} 
If ${\mathcal A}\subset E$ is tame then 
$A$ is tame and thus it has a mapping cylinder neighborhood 
$N(A)\subset F$. We could further assume that $F=N(A)$ by restricting to a smaller 
representative in $\mathfrak{E}_{n,q}^{PL}$ .  In fact, if 
$E'$ denotes $E\setminus (F\setminus N(A))\times S^{n-1})\subset E$, 
then $(E'\to S^{n-1}, {\mathcal A}) =(E\to S^{n-1}, {\mathcal A})\in \mathfrak{E}_{n,q}$. 

Further, we have $V=N(A)/A$.  As $N(A)\setminus A \cong \partial N(A)\times (0,1]$,  
it follows that $N(A)/A$ is homeomorphic to  the one point compactification of 
$\partial N(A)\times (0,1]$ and thus to the cone $C(\partial N(A))$. 
In particular $V$ is cone-like at $p$. 

Conversely, assume that $p$ is cone-like and let $N= C(L)$ be a conical neighborhood 
in $V$. Let then $E'=E\setminus (F\setminus N)\times S^{n-1})\subset E$, so that 
 $(E'\to S^{n-1}, {\mathcal A}) =(E\to S^{n-1}, {\mathcal A})\in \mathfrak{E}_{n,q}$.  
This way we  replaced $V$  in the new representative by a conical 
neighborhood $N$ and thus we can assume that $V=C(\partial V)$. Moreover,  
$\partial V= \partial F$ and thus $F/A$ is homeomorphic to $C(\partial F)$. 
We can therefore consider that $A$ is a spine of $F\setminus N(\partial F)$,  
where $N(\partial F)$ is a collar of $\partial F$ into $F$. In particular we can take 
$A$ to be tame, for any fiber $F$.   
This proves that  there is a representative of $f$ for which 
${\mathcal A}$ is tame.  
\end{proof}

The cone-like singularities were completely described in the 
second part of \cite{King}, where it is shown that any class 
in ${\mathfrak E}_{n,q}^{PL}$ is realized as $\delta(f)$ for some 
$f$ with a cone-like isolated singularity. 
The missing  part for completing the general 
classification (of not necessarily cone-like  
singularities) was the  characterization of the set of 
data in ${\mathfrak E}_{n,q}$ which could be realized as 
obstructions $\delta(f)$ for some $f$. 
We have an immediate but implicit description from above: 

\begin{proposition}
A class in $\mathfrak{E}_{n,q}$ can be realized by a map $f:M^{n+q}\to D^n$ 
with an isolated critical point  if and only if  the mapping cylinder 
$M(r)$ is a manifold (at $p$). 
\end{proposition}
\begin{proof}
Given the data $(E\to S^{n-1},  {\mathcal A}\subset E)$ we can recover both 
$M^{n+q}$ and the map $f$ as follows.  
The manifold $M^{n+q}$ (which is supposed to be an adapted neighborhood of 
$p$) is the mapping cylinder  $M(r)$ of the retraction $r$ i.e. 
$M(r)=E\times [0,1]\cup_{(x,1)\sim r(x)} V$. Moreover, the map 
$f:M^{n+q}\to D^n$ is given by the formula 
$f(x,t)=(f(x),t) $, where $D^n$ is identified with the mapping cylinder 
of the trivial map $S^{n-1}\to \{*\}$ to a point.  
This proves the claim. 
\end{proof}

\subsection{Contractibility of adapted neighborhoods and tameness}
The main result of this section gives a topological characterization  
of cone-like singularities in high dimensions:   
\begin{proposition}\label{contractible}
If the codimension is $q\leq 3$ then there exists a contractible 
adapted neighborhood $M^{n+q}$  
($n\geq 4$ if $q=3$) iff the singularity is cone-like. 
\end{proposition}
\begin{proof} 
One direction is clear. If the singularity is cone-like  the intersection $V\cap M$ 
with a small adapted neighborhood $M$ is $C(\partial V)$ and thus contractible. But 
$M$ deformation retracts onto $V\cap M$. 

In order to prove the reverse implication we need the following technical 
lemma:

\begin{lemma}\label{suite}
There exist open subsets 
$V_{s}\subset V$, $M_{s}\subset M$ and the sequence 
of concentric disks $D_s\subset D$ fulfilling the 
following conditions: 
\begin{enumerate}
\item $\overline{V_{s}}\setminus\{x\}$ is a manifold with compact boundary   
$\partial\overline{V_{s}}$. 
\item $V_{s}\subset V\setminus \partial V$.
\item Each component of $V\setminus\overline{V_{s}}$ meets $\partial V$. 
\item Each component of $V_{s}$ meets $x$. 
\item  $V_{s}=V\cap M_{s}$ and  the restriction 
$f:M_{s}\to D_{s}$ is proper with one critical point. 
\item  $\overline{M_{s}}$ is a compact manifold  
having the same components as $M_{s}$. 
\item  the diameter of $M_{s}$ tends to zero, when $s$ goes to infinity. 
\item We have $\overline{M_{s}} \subset M_{s+1}$. 
\item  The restriction $f:\overline{M}_{s}\cap f^{-1}(D_{s+1})
\setminus \overline{M}_{s+1}\to D_{s+1}$ is 
a fibration over the disk $D_{s+1}$ of radius.
\item We have  $\overline{D_{m}} \subset D_{m+1}$ and $\cap_{m=0}^{\infty} D_m=\{0\}$. 
\end{enumerate}
\end{lemma}
\begin{proof}
This is part of the content of  (\cite{CT3}, Lemma 3.3, p.951), where 
the lemma is stated for codimension 2. However, the same  
proof applies word-by-word without any codimension restriction for the 
statements above.  An alternative approach is given in the proof of Proposition 1 of \cite{King}, as an 
application of Siebenmann's Union Lemma 6.9 from \cite{Sie2}. 
\end{proof}

Let $F_{s}$ denote the generic fiber of 
$f:M_{s}\to D_{s}$ and $F$ be the generic fiber associated to $f:M\to D$.


The converse is given by: 
\begin{lemma}
If $\overline{M}$ is contractible and $n\geq 4$ then the 
 singular fiber $V\cap \overline{M}$ is homeomorphic to the 3-disk. 
\end{lemma}
\begin{proof}
Consider the sequences $V_s, M_s, D_s, F_s$ furnished by Lemma \ref{suite}  
 where $M_s$ are contractible. 
Lemma \ref{fiber} says that $F_s$ is a 3-disk, for all $s$. 

By hypothesis we have 
$\overline{M}_s\cap f^{-1}(D_{s+1}) \setminus \overline{M}_{s+1}\to 
D_{s+1}$ is a fibration and thus 
$(\overline{F}_s\setminus F_{s+1},\partial F_{s+1})$ is diffeomorphic to 
$(\overline{V}_s\setminus V_{s+1}, \partial \overline{V}_s)$.  
Therefore $\partial V_s$ is a 2-sphere and 
$\overline{V}_s\setminus V_{s+1}$ is a cylinder. 

The fiber $V\cap M\setminus \{x\}$ is the union 
$\cup_{s}(\overline{V_s}\setminus V_{s+1})$ i.e. 
the ascending union of cylinders. 
Since it has an exhaustion by compact sub-manifolds whose boundary 
are spheres it follows that  $V\cap M\setminus \{x\}$ is simply connected 
at infinity and thus it is diffeomorphic to $D^3\setminus \{0\}$. 
This implies that $V\cap M$  is homeomorphic to the 3-disk. 
\end{proof}
Proposition \ref{contractible} follows. 
\end{proof}


\subsection{Isolated singularities with trivial singular fibers after Hamstrom}
 The main result of this section concerns the singularities whose singular fiber is a disk. 

\begin{proposition}
If  the codimension is 3 and the singular fiber is homeomorphic to a 3-disk 
then $f$ is a topological fibration. 
\end{proposition}
\begin{proof}
Recall, after Hamstrom (\cite{Ham}) the following: 
\begin{definition}
An open proper 
mapping $f:X\to Y$ between metric spaces $X$ and $Y$ 
is homotopy $r$-regular if for each $x$ in $X$ and $\varepsilon>0$ there 
is a $\delta>0$ so that every mapping of a $s$-sphere $S^s$, 
with $s=0,1,\ldots,r$, into 
$D(x,\delta)\cap f^{-1}(y)$, for $y$ in $Y$, is homotopic to 0 in 
$D(x,\varepsilon)\cap f^{-1}(y)$. Here $D(x,\varepsilon)$ 
denotes the metric disk centered at $x$, consisting of points whose 
distances to $x$ are less than $\varepsilon$. 
\end{definition}

Assume now that we have an isolated singularity of codimension 3, as in the previous 
section. We know that  $\partial \overline{V}_s$ is a 2-sphere embedded in the disk $V$ 
and thus $V_s$ is homeomorphic to a 2-disk, by the Schoenflies Theorem.  
Since  the diameter of $\overline{V}_s$ 
goes to zero it follows that the singular fiber $V\cap M$ satisfies the 
homotopy 2-regularity condition at $0$. Thus $f$ is homotopy 2-regular.

It is proved in  (\cite{Ham}, Thm. 6.1)  that a homotopy 
2-regular map  $f$ of a complete metric space $X$ onto a 
finite dimensional space $Y$ such that each  fiber $f^{-1}(y)$ (for $y\in Y$)  
is homeomorphic to a given compact 3-manifold $Q$ with boundary without fake 3-disks
(i.e.  having the property that  each homotopy 3-cell in $M$ is homeomorphic to the 3-cell) 
is actually a  locally trivial fiber map. This proves the claim. 
\end{proof}

\begin{remark}
The  problem we faced above is the fact that the singularity might not be 
cone-like at $x$.  If the singularity were tame in the sense developed 
by Kauffman and Neumann  in \cite{KN}
(for instance $f$ is locally  topologically 
equivalent to a real polynomial mapping) then we could take 
for $M^{n+q}$ a ball and  the arguments in 
\cite{Milnor} would immediately show that the local Milnor fiber $F^3$ is  
contractible (when dimension is $n\geq 4$) or simply connected if 
$n=3$. Thus $F^3$ is a disk and thus the singularity is trivial, as proved 
in \cite{CL}. It is actually proved in (\cite{King}, p.396) that 
cone-like singularities have adapted neighborhoods which are 
homeomorphic to balls at least if the dimension 
$n+q$ is not 4 or 5.  

Now, if the singularity is not cone-like,  we can always assume that the 
singular fiber $V$ is transverse to all small enough spheres 
$S_{\varepsilon}$ centered at $x$, but for an infinite discrete set 
of radii $\varepsilon_n$ accumulating to $0$. However, we have no 
local Milnor fibration and no local cone structure around a singularity. 
We have instead approximations of this fibration, as provided 
by Lemma \ref{suite}.
\end{remark}

\subsection{Removable singularities}

\begin{definition}\label{remov}
A smooth map $f:M\to D$ with a single critical point has a 
{\em removable} isolated singularity 
at $p$ if there exists a  smooth manifold $X$, a smooth map 
$g:X\to D$ without critical points 
and a homeomorphism $h:M\to X$ which is a diffeomorphism in a 
neighborhood of the boundary 
such that  outside of an open neighborhood  of the critical point 
we have $f=g\circ h$.  
The singularity is  called {\em strongly removable} if one can  
obtain $X$ by excising a  smooth disk  out of  $M$ 
and gluing it back.  
\end{definition}

A typical example of a map with a removable isolated singularity 
is that obtained from a map $g$ without critical points by composing with 
a {\em near-homeomorphism}. 
Recall that near-homeomorphisms are limits of homeomorphisms 
in the compact-open topology and they were characterized by 
Siebenmann (\cite{Sie}) as being the cell-like maps between manifolds.

\begin{proposition}
If the codimension is at most 2  then isolated singularities are removable. 
\end{proposition}
\begin{proof}
Assume that we have an orientable fibration $E\to S^{n-1}$ having the 
fiber $F^2$ of dimension at most 2. 

If the codimension is 0 then  either $n=2$ and 
$E$ is a finite covering and thus $A\subset F$ is obviously tame 
since it is a finite set of points or else $n\geq 3$ and thus 
$F$ is a point and $A$ is empty. 

If the codimension is 1 then the fiber is a bunch of intervals. As the boundary 
fibration must be trivial it follows that $E$ is a product.

If the codimension is 2 then $F$ is an orientable surface with boundary. 
A homeomorphism which is identity on the boundary should preserve 
the orientation. Further, a celebrated result due to Earle, Eells and Schatz 
tells us that all connected components of 
${\rm Homeo}^+(F,\partial F)$  are contractible 
(see \cite{EE, ES}), if $F$ is neither a disk nor an annulus. 
In particular, such a  surface fibration over 
$S^{n-1}$ (with $n\geq 2$)  must be trivial. 
This implies that  either $E$ is trivial or the 
fiber is a disk or an annulus. 

Moreover, if the surface is either a disk or an annulus 
then ${\rm Homeo}^+(F,\partial F)$ 
has the homotopy type of a circle.  In particular 
a  disk or annulus fibration over 
$S^{n-1}$ is trivial whenever $n\geq 3$.

If the fibration $E$ is trivial then the mapping cylinder $M(r)$ of any 
map $r:E\to F/A$ can be  globally described as a quotient.  
In fact, $M(r)$ is given by the quotient 
$F\times D^n/A\times \{0\}$.  This is the simplest example 
of a decomposition space obtained from a upper semi-continuous 
decomposition, which 
in our case  has only one non-degenerate element, namely $A\times \{0\}$. 

The theory of decompositions of manifolds is largely described in 
the Daverman monograph \cite{Dav}. 
According to (\cite{Dav}, Exercise 7, p.41)   
$F\times D^n/A\times \{0\}$ is a manifold 
iff $A\times \{0\}\subset F\times D^n$ is cellular. 
Recall that a subset $A$ of an $n$-manifold 
is called {\em cellular} if there exists a nested  sequence of 
$n$-cells  $B_i$ in that manifold  such that 
$B_j\subset {\rm int}(B_{j-1})$ for any $j$ and 
$A=\cap_{j=1}^{\infty} B_j$. 

\begin{remark}
Cellular sets are compact and connected but not 
necessarily locally connected, as 
the graph of $\sin\left(\frac{1}{x}\right)$ is such a 
cellular subset of $\R^2$. 
\end{remark}

Moreover,   it is known that  the decomposition 
consisting of a single cellular set is shrinkable  
(see \cite{Dav}, Corollary 2.2A, p.36)  and thus the quotient of a 
manifold by a cellular subset 
is again a manifold 
(see \cite{Dav}, Theorem 2.2, p.23) homeomorphic to the former 
manifold. 

Alternatively,  a cellular subset $A$ in a manifold is {\em cell-like} 
(see \cite{Dav}, Corollary 3.2.B, p.120) 
namely,  for each neighborhood 
$U$ of $A$ in that manifold one can contract  $A$ to a point in 
$U$.  Compact ANR are cell-like if and only if they are contractible.  
Cell-like maps are those continuous maps whose point inverses 
are cell-like sets. 
Siebenmann characterized in \cite{Sie} the cell-like maps between manifolds 
of the same dimension as being precisely those maps which are 
near-homeomorphisms. 
Therefore the projection  map $F\times D^k\to F\times D^k/A$, which 
is a cell-like map, is a near-homeomorphism.  

In particular we can replace the map 
$M\to D^n$ by composing with the near-homeomorphism 
$F\times D^n\to M$ by the projection $F\times D^n\to D^n$ 
which has no critical points and has the same boundary restrictions. 
In other words the singularity of the map $f$ is removable.

Otherwise the manifold $F$ is a 2-disk or an annulus. 
Recall that we actually have a nested sequence of fibers 
$F_s\subset F_{s-1}\subset F$, 
whose intersection is the vanishing compactum $A\subset F$ in the fiber. 

If all but finitely many $F_s$ are disks then 
$A$ is cellular in $F$.  From \cite{Dav}
$A$ is shrinkable and thus $V=F/A$ is homeomorphic to a disk. 
We can apply 
Hamstrom's Theorem from above to obtain that 
$f$ has no topological critical points.   

If all but finitely many $F_s$ are cylinders and $n=2$ then   
$A$ is a circle in $F$.  
In fact the boundary of the cylinders $\partial F_s$ should 
be parallel and thus they consist in arbitrarily thin adapted neighborhoods 
of the  intersection circle   $A=\cap_{s=0}^{\infty} F_s$. Thus $A$  is tame and 
hence the singularity is cone-like. 

If $F$ is non-orientable then the group of homeomorphisms is $SO(3)$ (for $S^2, \R\mathbb P^2$)
or else a circle (for the Klein bottle, disk, annulus or the Mobius band). 
The argument  above for the annulus extends for the Mobius band as well. 
\end{proof}

\begin{remark}
We can actually show that all codimension 2 singularities are cone-like, as 
claimed in Proposition \ref{codim2}. We give here an alternative proof below. 
\end{remark}
\begin{proof}
It suffices to show that $A\subset F$ is tame by constructing a 
regular neighborhood of it. 
Recall that $A=\cap_{s=0}^{\infty}F_s$, where 
$F_{s+1}\subset {\rm int}(F_s)$, for all $s$. 
Moreover, we can show that the number of components of $\partial F_s$ 
is uniformly bounded. 
It is clear that genus of $F_s$ is decreasing. 
Thus a descending sequence of surfaces 
like above has the homeomorphism type 
eventually constant. Further by finiteness and  the Gruschko 
Theorem  it follows that the surfaces have eventually parallel 
boundaries. This implies that 
we have a mapping cylinder neighborhood for $A$. 
\end{proof}

\begin{proposition}\label{codimthree}
\begin{enumerate}
\item In codimension 3 there exist uncountably many non-equivalent 
germs of continuous maps with an isolated singularity. 
\item 
Any codimension 3 smooth map with an  isolated singularity $M\to D^n$ 
 for $n\geq 5$ is  either  strongly removable or else an isolated point in 
the fiber. 
\end{enumerate}
\end{proposition}
\begin{proof}
The first point is to observe that there exist $A\subset F^3$ which are wild. 
For instance we can take any  Artin-Fox wild arc (see \cite{ArF}) in $D^3$ or 
a Whitehead continuum such that its complement in 
$S^3$ is a Whitehead manifold. 
It is known (see \cite{Br}) that there exist uncountably many 
(one-holed)  Whitehead manifolds $W^3$ which are pairwise non-homeomorphic. 
According to the King classification of germs it follows that 
maps associated to the data 
$(E=D^3\times S^{n-1}, {\mathcal A}=(D^3-W^3)\times S^{n-1})$ 
are not right equivalent if the  respective $W^3$ are not homeomorphic. 
On the other hand the mapping cylinder is 
$M=D^3\times D^n/A\times \{0\}$.  

Let  $A$ be an arbitrary  contractible closed subset of $D^3$,  
for instance a wild arc.  It follows that 
the decomposition  given by the singleton $A\times\{0\}$ of $D^3\times D^n$ 
 is shrinkable, as soon as $n\geq 2$ (see \cite{Dav}, Corollary 8A, p.196). 
Therefore  $M$ is homeomorphic to $D^3\times D^n$ and thus it is a manifold.  
In particular, we obtain uncountably many germs of continuous maps 
$R^{n+3}\to \R^n$ 
which are pairwise  
non-equivalent by right composition with a homeomorphism.

Let us prove now the second claim, namely that all codimension 3 
singularities are arising as claimed.  We suppose henceforth 
that the singular fiber is not reduced to one point. 

If the fiber $F^3$ is irreducible and the boundary is non-empty then 
either $F^3$ is a disk or else 
${\rm Homeo}^+(F,\partial F)$ has contractible components, according to 
Hatcher (see \cite{H,Hat,HaM}).  Therefore, either $F$ is a 3-disk 
or else $F^3$-fibrations over $S^{n-1}$ are trivial, as soon 
as $n\geq 3$. 

Let us analyze the case where $F$ is not irreducible. Recall that we have 
a sequence of nested  manifold neighborhoods $F_{s+1}\subset F_s\subset F$. 
If there exists at least one $F_s$ which is irreducible then we can apply the argument from above. 

\begin{lemma}\label{irred}
If $n\geq 5$ then there exists some representative of the germ $f$ 
of isolated singularity for which 
some associated fiber $F^3$ is irreducible. 
\end{lemma}
\begin{proof}
Suppose that all $F_s$ are reducible. 
According to the Sphere Theorem 
(see \cite{Hempel,Mil2}) we can detect the essential 2-spheres 
as follows: $F_s^3$ is reducible  if and only if $\pi_2(F_s)\neq 0$, modulo the 
Poincar\'e Conjecture. 

Therefore, by the Haken-Kneser 
finiteness Theorem there exist only finitely many  
embedded  essential and  pairwise non-parallel 2-spheres 
$S_{s,t}\subset F_s$. For fixed $s$ the spheres $S_{s,t}$ are disjoint. 
 Further, if for some pair $s_1 < s_2$ we have 
$S_{s_1,t}\cap S_{s_2,v}\neq \emptyset$, then we can replace 
$S_{s_2,v}$ by a number of spheres obtained by pushing them along the 
intersections. 
In fact the intersection of two spheres $S_1$ and $S_2$ can be made transverse 
by a small isotopy and $S_1\setminus S_1\cap S_2$ is a genus zero surface with 
boundary a number of circles. Each such circle, when viewed in $S_2$, bounds a disk in $S_2$.
We start from  a circle bounding an innermost such disk $\delta$ and then consider  two parallel copies 
of $\delta$ corresponding to the  parallel boundary circles of a cylinder neighborhood $U$ of the 
circle in $S_1$. 
We remove then the cylinder $U$ from the sphere $S_1$ and glue back the two copies 
of $\delta$. This replaces $S_1$ by two spheres whose product as class in $\pi_2(M)$ 
is the class of $S_1$. Moreover the intersection 
between $S_1$ and $S_2$ countains one circle less than before. 
We can continue this procedure until $S_2$ gets disjoint from the 
spheres we create.  

There are finitely many such spheres, and among them there exists at least one which is 
non-trivial in $\pi_2(F_{s_2})$ since the product of all the 
spheres obtained this 
way is the class of the original sphere $S_{s_2,v}$. 
By using a double recurrence on $s_2$ and $s_1$ 
we obtain that the spheres $S_{s_1,t}$ are 
pairwise disjoint for all $s_1$ and $t$. 

Recall now that we have a 
fibration of $M_s\setminus V_s\to D^{n}\setminus\{0\}$ with fiber $F_s$. 
Its homotopy  exact sequence reads: 
\[ \pi_{i+1}(S^{n-1})\to \pi_i(F_s) \to \pi_{i}(M_s\setminus V_s)\to \pi_i(S^{n-1})\to \]
From Lemma \ref{connect} and the proof of Lemma \ref{fiber} 
the inclusion $F_s\hookrightarrow M_s$ induces 
$\pi_i(F_s)\cong \pi_i(M_s)$, for $i\leq 2$, as soon as $n\geq 5$. 
Since $M_s$ deformation retracts onto $V_s$ the composition 
$F_s\hookrightarrow M_s\to V_s$ induces an isomorphism:  
\[ \pi_i(F_s)\cong \pi_i(F_s/A), \;\; {\rm for }\,\; i\leq 2\]
We will show below that this is contradicted for $i=2$ when all 
$F_s$ are reducible.

Using again the Haken-Kneser finiteness Theorem 
there exists some $s$ such that one of the following holds: 

\begin{enumerate}
\item either there exists some infinite 
sequence  of spheres $S_{k_i,v_i}$ which are 
parallel to some $S_{s,t}$ within $F_s$;

Let now $U_i$ be the annuli bounding 
$S_{k_i,v_i}\cup S_{k_{i+1},v_{i+1}}$ in $F_s$  
and $O_i$ be their image in $F_s/A$, where $S_{k_0,v_0}$ denotes $S_{s,t}$. 
Set $O_{\infty}=\cup_{i=0}^{\infty} O_i \cup A\subset F_s/A$  be the union 
of the annuli and the singular point $A$. As $\cup_{i=0}^{\infty} U_i$ is a 
punctured 3-disk and $S_{k_i,v_i}$ accumulate on $A$ 
it follows that $O_{\infty}$ is  the  continuous image of  the 
compactified 3-disk into $F_s/A$. 
This implies that the image of the class of $S_{s,t}$ under the map  
$\pi_2(F_s)\to \pi_2(F_s/A)$ vanishes. But the class of 
$S_{s,t}$  in $F_s$ is non-zero since this is an essential sphere.  
This contradicts the isomorphism above. 

\item or else the 2-spheres $S_{k,v}$, for $k\geq s+p$ 
are null-homotopic in $F_s$, for some $p$. 

Then $S_{k,v}$ is essential in $F_{k}\subset F_s$ 
but null-homotopic in $F_s$. By the Sphere Theorem, possibly changing 
the collection of essential spheres, we can assume that each 
$S_{k,v}$ bounds an embedded 3-ball $B_{k,v}$ in $F_s$. 
Let $M_k$, $E_k\to S^{n-1}$ be the adapted neighborhood and respectively 
fibration of the germ $f$, corresponding to the fiber $F_k$. 
Then $\overline{M_k}=M_k\cup_{v} B_{k,v}\times D^{n}\subset M_s$, because 
$M_s\setminus M_k$ is a trivial fibration over the disk $D^n$. 
Let $k\geq s+p$  and consider the restriction of the germ 
$f:M_s\to D^n$ to $\overline{M_k}$. The restriction of the fibration 
$E_s$ to $E_s\setminus E_k$ is a trivial fibration over $S^{n-1}$ 
with fiber $F_s\setminus F_k$, because 
it extends over $D^n$. Thus $\overline{E_k}=
E_s\cap \overline{M_k}$ is a sub-fibration of $E_s$ over 
$S^{n-1}$ whose fiber is $\overline{F_k}=F_k\cup_{v} B_{k,v}$. 
The fibration $\overline{E_k}\to S^{n-1}$ (with the same vanishing compactum) 
defines therefore the same germ $f$.  

Now the 2-spheres $S_{k,v}$ are pairwise 
disjoint and $\partial F_k\subset \cup_v B_{k,v}$. 
Therefore every essential 2-sphere $S$ contained in the interior of   
$F_k\setminus \cup_v B_{k,v}$ is parallel to one of the boundary 
spheres $S_{k,v}$. In fact such a 2-sphere should be homotopically 
trivial in $F_s$ and thus, by the Sphere Theorem we can 
suppose  that $S$ bounds a ball $B$ in $F_s$. Thus either $B$ is disjoint 
from  the balls $B_{k,v}$ or there exists some $v$ so that 
$B_{k,v}\subset {\rm int}(B)$. In the first case $S$ is trivial in 
$F_k\setminus \cup_v B_{k,v}$ and in the second case we could use Brown's 
solution to the Schoenflies Conjecture to obtain that 
$S\cup S_{k,v}$ bounds an annulus in $F_k$ and hence in 
$F_k\setminus \cup_v B_{k,v}$.  
Eventually $\overline{F_k}$ is irreducible since all essential 
spheres of $F_k\setminus \cup_v B_{k,v}$ were capped off by 3-balls. 
Thus, in this case there exists a representative of the germ $f$ whose 
fiber ${\overline F_k}$ is irreducible. 
\end{enumerate}
\end{proof}

Now we can conclude as above: either $F_s$ are disks for  all large enough $s$ 
or else the  $F_s$-fibrations are trivial. If all $F_s$ are 3-disks then 
$A$ is cellular in $F$ and thus 
the map $F\to F/A=V$ is a near-homeomorphism and thus the singular fiber is 
homeomorphic to a disk.  
Then Hamstrom's Theorem from the previous section implies that 
the singularity is removable. 

If some $F_s$ fibration is trivial it follows that $M(r)$ is $F_s\times D^n/A\times \{0\}$. 
From \cite{Dav} $M(r)$ is a manifold only if $A\times \{0\}$ is cellular in $F_s\times D^n$ 
and in this case the projection $F_s\times D^n\to F_s\times D^n/A\times \{0\}$ is a 
near-homeomorphism. 

Moreover, in both cases we are able to be more precise. As 
$A\times \{0\}$ is cellular 
it admits a sequence of arbitrarily small disk neighborhoods in $F\times D^n$. 
Furthermore one can choose then a  PL and thus a smooth  neighborhood $B$ of it.
Then $B/A\times \{0\}$ is also homeomorphic to a disk since it is 
a quotient by a cellular subset and the restriction of the 
map $B\to B/A\times \{0\}$ to a 
collar of the boundary is the identity. 
Thus the boundary of $B/A\times \{0\}$ is a smooth sphere embedded   
in the smooth manifold $M$ and hence $B/A\times \{0\}$ is a smooth manifold homeomorphic 
to a disk. As $n\geq 2$ Smale's Theorem implies that 
$B/A\times \{0\}$ is diffeomorphic to a disk. 

Let $q:B\to B/A\times \{0\}$ be a homeomorphism 
(recall that the projection is a near-homeomorphism), 
which can be assumed to be smooth on the 
boundary. It follows that 
we can obtain the manifold  $F\times D^n$  from $M$ by excising a $(n+3)$-disk 
corresponding to $B/A\times \{0\}$ and gluing it back by twisting the gluing map 
by the homeomorphism $q$; this is the same as gluing back $B$.  

Therefore, with the notations from Definition \ref{remov}, one  
puts $X= F\times D^n$ to  
find that the singularity is strongly removable, as claimed. 
\end{proof}

\begin{remark}
It seems that isolated singularities in dimensions 
$(6,3)$ and $(7,4)$ are  either cone-like or else removable. 
From above it suffices to analyze the case 
when $F_s$ are reducible. If $n=4$ and we knew that there are adapted 
neighborhoods which are contractible, then the claim would follow. 
\end{remark}

\begin{remark}
Examples in the next section will show that the statement of 
Lemma \ref{irred} cannot be extended to $n=3$, since there exist fibered links 
whose fibers are not irreducible. 
Moreover, our proof of Proposition \ref{codimthree} 
could not work for $n=3$ since the fibered links 
fibrations might also be non-trivial. 
\end{remark}

\subsection{Cone-like isolated singularities in dimensions $(6,3)$}

There exists a general recipe for constructing fibered links 
$\sqcup_{d}S^{n-1}$ into $S^{2n-1}$. We will restrict to the high 
dimensional case $n\geq 3$ below. According to 
Haefliger (see \cite{Hae2a,Hae2b}) embeddings of 
$S^{n-1}$ into $S^{2n-1}$ are trivial up to isotopy, when 
$n\geq 3$. Thus to each $S^{n-1}_1\subset S^{2n-1}$ one can associate 
a {\em Hopf dual} $S^{n-1}_2$ which links $S^{n-1}_1$ once. 
In fact $S^{n-1}_2$ is the core of the complement of a 
regular neighborhood of $S^{n-1}_1$. 

Consider then $S^{n-1}_2, \ldots S^{n-1}_{d}$ be a set of pairwise 
disjoint Hopf duals  to $S^{n-1}_1$. Then, by Haefliger classification Theorem 
(see \cite{Hae2a,Hae2b})  the link $L=\sqcup_{j=1}^d S^{n-1}_j$ is 
uniquely determined, up to isotopy,  
by its linking matrix $A_L$. Since the first column and row are 
made of $\pm 1$, the most important information is the linking matrix 
$A_{L^*}$ of the sub-link  $L^*=\sqcup_{j=2}^d S^{n-1}_j$. 
Observe also that links were oriented for the purpose of 
computing linking numbers, although the fibering property 
concerns unoriented links. 

By convenience the diagonal of the matrix $A_L$ will have trivial entries 
$0$ as the links are not framed. We also set $L_A$ for the link $L$ 
with $A_{L^{*}}=A$.
The linking numbers $lk(a,b)$ of the $(n-1)$-spheres $a,b$ 
in $S^{2n-1}$ satisfy $lk(a,b)=(-1)^{n}lk(b,a)$.

\begin{proposition}
\begin{enumerate}
\item For every  choice of an integral 
$(-1)^n$-symmetric $(d-1)\times (d-1)$ matrix $A$ with trivial diagonal 
the link $L_A$ is fibered. 
\item If $n=3$ then every fibered link with simply connected fiber 
is isotopic to some $L_A$.    
\end{enumerate}
\end{proposition}
\begin{proof}
The complement of a regular neighborhood
of an $S^{n-1}$ in $S^{2n-1}$ is $S^{n-1}\times D^n$, which fibers 
trivially over $S^{n-1}$. 

Let us consider now sub-fibrations 
of $S^{n-1}\times D^n$ for which the first projection restricts 
to a locally  trivial fibration over $S^{n-1}$ with fiber 
$D^{n}\setminus \sqcup_{i=1}^{d-1}D^n$.  
Such locally trivial fibrations $\eta$
are determined up to  isotopy among fibrations  
by the homotopy  class of the  classifying 
map $c_{\eta}:S^{n-1}\to {\mathbb F}_{d-1}(D^n)$, 
where ${\mathbb F}_{d-1}(D^n)$ denotes 
the configuration space of $(d-1)$ disjoint disks of 
equal (very small) radii in $D^n$. The value 
$c_{\eta}(x)$  is the class in ${\mathbb F}_{d-1}(D^n)$ of the fiber  
$\eta(x)\subset D^n$. Moreover ${\mathbb F}_{d-1}(D^n)$ 
has the homotopy type of the configuration space 
${\mathbb F}_{d}(\R^n)$ of $d$ points on $\R^n$. 
 
The map ${\mathbb F}_{d}(\R^n)\to {\mathbb F}_{d-1}(\R^n)$ 
which forgets the first point of the $d$-tuple is a fibration 
whose fiber ${\mathbb F}_{1,d-1}(\R^n)$ has the homotopy type of a wedge 
$\vee_{d-1}S^{n-1}$ of $(d-1)$ spheres of dimension $(n-1)$ (see \cite{FH}).  
In particular,  ${\mathbb F}_{d-1}(\R^n)$ is $(n-2)$-connected.

To every map ${\bf f}:S^{n-1}\to {\mathbb F}_{d}(\R^n)$ one associates 
$d$ disjoint embeddings $f_j:S^{n-1}\to S^{n-1}\times D^n$ 
by the formula $f_j(x)=(x,{\bf f}(x)_j)$, where 
${\bf f}(x)_j$ is the $j$-th component of the $d$-tuple 
${\bf f}(x)\in {\mathbb F}_{d}(\R^n)$, and one identifies 
${\rm int}(D^n)$ to $\R^n$.

\begin{lemma}
We have an isomorphism 
$\pi_{n-1}({\mathbb F}_{d-1}(\R^n))\to\Z^{(d-2)(d-1)/2}$ which 
associates to the class of the map 
${\bf f}:S^{n-1}\to {\mathbb F}_{d-1}(\R^n)$ the $(-1)^n$-symmetric 
matrix whose upper triangular  entries are given by 
$(lk(f_i(S^{n-1}),f_j(S^{n-1})))_{1\leq i < j\leq d-1}$.  
\end{lemma}
\begin{proof}
It is known that $\pi_{n-1}({\mathbb F}_{d-1}(\R^n))$ is the direct sum 
$\pi_{n-1}({\mathbb F}_{d-1}(\R^n)))\oplus \pi_{n-1}({\mathbb F}_{d-2}(\R^n))$. 
By functoriality it suffices then to consider  
only the case when $d=2$, where the proof is elementary. 
In fact, let $\pi$ denote 
the projection  ${\mathbb F}_{2}(\R^n)\to S^{n-1}$ given by  
$\pi(x,y)=\frac{x-y}{|x-y|}$. Then $lk(f_1(S^{n-1}),f_2(S^{n-1}))$ 
equals the degree of $\pi({\bf f})$, as an element of $\pi_{n-1}(S^{n-1})$. 
\end{proof}

It follows that for every $(-1)^n$-symmetric matrix 
$A$ with trivial diagonal there exists some sub-fibration 
$\eta$ with fiber $D^{n}\setminus \sqcup_{i=1}^{d-1}D^n$ 
whose associated link is $L_A$. Thus $L_A$ are fibered links. 

Consider next $(m,n)=(6,3)$. Then the fiber $F^3$ of a fibered link is 
simply connected and hence a disk-with-holes 
$D^{3}\setminus \sqcup_{i=1}^{d-1}D^3$. The associated fibration is 
then a sub-fibration of a $D^3$-fibration over $S^2$, which is 
trivial (by the Hatcher solution to the Smale Conjecture in the smooth 
category and by the Alexander trick in the PL and topological categories). 
Therefore the fibered link arises as in the previous construction. 
Observe that the singular fiber is 
the cone over the union of $d$  disjoint 2-spheres, namely a 
wedge of $d$ disks of dimension $3$.  
\end{proof}

\begin{remark} 
The proof  from above shows that the isolated singularity so obtained  
is locally topologically a fibration, if $d=1$ and $n=3$. 
\end{remark}

Once we know all fibered links in dimension $(6,3)$ we can construct 
examples of  pairs of manifolds with finite $\varphi$ by using 
a Lego construction. Let $\Gamma$ be a (decorated) graph whose vertices have 
valence at least 3. Each vertex $v$ of $\Gamma$ is decorated by some  
  symmetric integral $(d-1)\times (d-1)$ matrix $A(v)$. To every vertex 
$v$ there is associated a fibration 
$S^5-N(L_{A(v)})\to S^2$ which extends to a smooth map with 
one critical point $f_v:D_v^6\to D^3$. We glue together the disks 
$D_v$ using the pattern of the graph $\Gamma$ by identifying 
one component of $\partial N(L_{A(v)})$ to one component 
of $\partial N(L_{A(w)})$ if $v$ and $w$ are adjacent in $\Gamma$. 
The identification has to respect the trivializations  
$\partial N(L_{A(v)})\to D^3$ and hence one can take them to be 
the same as in the double construction. 
We obtain then a manifold with boundary $X(\Gamma, A(v)_{v\in \Gamma})$ 
endowed with a  proper smooth map $f_{\Gamma}:X(\Gamma, A(v)_{v\in \Gamma})\to D^3$ with $n$ singular points 
($n$ being the number of vertices of $\Gamma$) inside one fiber. 
The generic fiber is  $\#_{g}S^1\times S^2$, where $g$ is the rank of 
$H^1(\Gamma)$. If orientation reversing  
gluing homeomorphisms are allowed then  one 
could also obtain  non-orientable factors 
homeomorphic to the twisted $S^2$-fibration over the circle. 
The restriction of $f_{\Gamma}$ to the boundary is a 
$\#_{g}S^1\times S^2$-fibration over $S^2$. Let 
$\Gamma_1,\Gamma_2,\ldots,\Gamma_p$ be a set of graphs associated to 
a family of cobounding fibrations, 
namely such that there exists a fibration over 
$D^3\setminus \sqcup_{i=1}^{p-1}D^3$ extending the boundary 
fibrations restrictions of  $f_{\Gamma_i}$, $1\leq i\leq p$. 
Then we can glue together $f_{\Gamma_j}$  to obtain some smooth map 
with finitely many critical points into $S^3$. 
In particular, we can realize the double 
of $f_{\Gamma}$ by gluing together 
$f_{\Gamma}$ and its mirror image. 

\begin{remark}
All 6-manifolds $M^6$ admitting a smooth map $M^6\to S^3$ 
with finitely many cone-like singularities should  
arise by a similar construction. However, at present we don't know whether 
there are any non-trivial examples, in particular we are unable  
to find whether $\varphi(M^6,S^3)$ equals the total 
number of vertices of the graphs $\Gamma_j$. Notice that 
a similar construction produces the examples in 
dimensions $(4,3)$ and $(8,5)$, but in the later case the value 
of $\varphi(M,N)$ can be expressed 
in terms of algebraic topology invariants of the manifolds $M,N$.
\end{remark}

\subsection{Essential singularities in codimension 3 and proof of Theorem \ref{local}}
Let $M^{n+3}\to N^{n}$ be a smooth map with finitely many critical points. 
We have to prove that  
$M^{n+3}$ is diffeomorphic to $\Sigma^{n+3}\# Q^{n+3}$ where:  
\begin{enumerate}
\item $\Sigma^{n+3}$ is a homotopy sphere; 
\item  there is a smooth  (induced) map $Q^{n+3}\to N^n$ having only finitely many critical points and 
fulfilling:
\begin{enumerate}
\item if $n\geq 6$ or $n=4$ and the singularities were cone-like 
then there are no critical points; 
\item  if $n=5$ then around  each critical point the map is smoothly equivalent to the Kuiper map 
$k_{4,1}$ i.e. the cone over the  corresponding Hopf  fibration. 
\end{enumerate}
\end{enumerate}

Around each critical point of $f:M\to N$ we can find adapted neighborhoods 
$M_i\subset M$ and disks $D_i\subset N$ centered at the critical points so that 
$f:M_i\to D_i$ is a proper smooth map with precisely one critical point, 
for $i\in\{1,2,\ldots,r\}$.

A Theorem of Timourian (see \cite{Tim}) states that a critical point 
which is an isolated 
point of its fiber is such that locally the map is  topologically 
equivalent to the cone over a Hopf fibration (and thus $n=5$). 
Thus the $M_i$, $1\leq i\leq s$,  
corresponding to such critical points are homeomorphic to 
a disk $D^{n+3}$ and thus diffeomorphic to $D^{n+3}$, 
by Smale's Theorem, for $n\geq 2$.

From the strong removability of the  remaining critical points  we can find maps 
$f_i:M_i'\to D_i$, where $M_i$ and $M_i'$ are homeomorphic, for all $i$ 
with $s+1\leq i\leq r$. Moreover 
$M_i\setminus U_i$ is identified with $M_i'\setminus U_i'$ for some open disks  
$U_i$ and $U_i'$ and then $f_i$ and $f|_{M_i}$  coincide when 
restricted to the subset $M_i\setminus U_i$. 
Therefore we can glue together the various mappings  
$f_i$ to obtain a smooth map $f': M'\to N$, where 
$M'=  \cup_{i\geq s+1} U_i' \cup \left(M\setminus \cup_{i\geq s+1} U_i\right) $. 
This amounts to excise a number of disjoint disks and glue them back differently. 
The result is then the connected sum of $M$ with a homotopy sphere. 

We excise now each  $M_i\subset M'$  for all $i\leq s$ and glue it back by using 
the restriction of the homeomorphism which identifies $f|_{M_i}$ with the cone over 
the Hopf fibration. This amounts to change $M$ by making one more connected sum with a 
homotopy sphere, since these $M_i$ (with $i\leq s$) are disks. 
Let then $M''$ be the resulting manifold. Then the induced map $M''\to N$ is smooth and 
has no other critical points but the $s$ points where it is smoothly equivalent to  the 
Kuiper map. Moreover, $M''$ is homeomorphic to $M$ and obtained from $M$ by 
a connected sum with a homotopy sphere.

\section{Proof of the structure Theorem \ref{complete}}
\subsection{Spherical blocks and Montgomery-Samelson fibrations}
Here and henceforth the dimensions $m,n$ will be of the form 
$m=2k, n=k+1$, where $k\in\{2,4\}$, unless the opposite is explicitly stated. 
The subject of this section is the description of the structure of 
spherical blocks, namely of manifolds $M^{2k}$ admitting a smooth 
map with finitely many critical points into $S^{k+1}$.   
Specifically, we have first: 

\begin{proposition}\label{char}
Assume that $M^{2k}$ is a compact orientable manifold.  
\begin{enumerate}
\item If $(2k,k+1)=(4,3)$ and
$\varphi(M^4,S^{3})$ is finite non-zero then
$M^4$ is either homeomorphic to $\#_r S^{2}\times S^{2}$, for some $r\geq 0$ or 
else homeomorphic to a fibration over $S^3$. 
\item  If $(2k,k+1)=(8,5)$ and
$\varphi(M^8,S^{5})$ is finite non-zero then $M^8$ 
is diffeomorphic either to
\begin{enumerate}
\item 
$\Sigma^{8}\#_r S^{4}\times S^{4}$ for some $r\geq 0$, with $\Sigma^8$ a
homotopy sphere, or to
\item $\Sigma^8\# N^8$, where $N\to S^5$ is a fibration,
$\Sigma^8$ is an exotic sphere (actually the generator of  the
group of homotopy 8-spheres $\Gamma_8=\Z/2\Z$), while
$\Sigma^8\# N^8$ is not a fibration over $S^5$.
\end{enumerate}
\end{enumerate}
\end{proposition}

Recall now the following definition from \cite{MonSam}:
\begin{definition}
A  Montgomery-Samelson
fibration $f: M^{m}\rightarrow N^{n}$
with singular set $A\subset M$ is a smooth map whose 
restriction to $M-A$
is a locally trivial fiber bundle while  the
restriction to $A$ is a homeomorphism.
\end{definition}

As a consequence of Proposition \ref{dim} and Theorem 
\ref{local} it follows
that, if $m-n\leq 3$ and $(m,n)\not\in\{(4,2), (5,2), (6,3)\}$ 
then  a smooth map $f:M^m\to N^{n}$ with finitely many critical 
points is a Montgomery-Samelson fibration with finite set $A$.
We used the result of Timourian (\cite{Tim}), stating 
that  $x$ is not an isolated point of $f^{-1}(f(x))$ unless  
$f$ is locally topologically equivalent to the suspension of a 
Hopf fibration.

Moreover, Montgomery-Samelson fibrations  of closed
orientable manifolds $M^m$ over spheres
were completely characterized topologically by
Antonelli, Church, Timourian and Conner  (see \cite{Ant3,CT1,CT2,Con})
as follows:

\begin{proposition}\label{act}
Let $M^m\to N^{n}$ be a Montgomery-Samelson fibration  with finite
non-empty singular set $A$ and $n\geq 2$.
\begin{enumerate}
\item  Then $(m,n)\in\{(2,2), (4,3), (8,5), (16,9)\}$ and, if the 
codimension is positive then the fiber is a sphere. In particular 
$m=2k, n=k+1$, for $k\in \{1,2,4,8\}$.
\item  If $M^m$ is simply connected then $|A|=\chi
(M)=\beta_{k}(M) +2$.
\item  When $N=S^{k+1}$ then $M^{2k}$ is homeomorphic to
$\#_{r}(S^{k}\times S^{k})$, where
$r=\frac{1}{2}\beta_{k}(M)$.
\end{enumerate}
\end{proposition}
\begin{proof}
See \cite{Ant1,Ant3,Tim}. When $N=S^{2}$ it was only shown in \cite{Ant3}
that  $M^4$ has the oriented homotopy
type of $\#_{r}(S^{2}\times S^{2})$, but the classification Theorem
for topological 4-manifolds due to Freedman 
permits to conclude. We used here the Poincar\'e Conjecture
for the 3-dimensional fiber.
\end{proof}

{\em Proof of Proposition \ref{char}}.
Consider now that $(m,n)\in \{(4,3), (8,5)\}$.
If $\varphi(M^{2k},S^{k+1})$
is finite non-zero then $M^{2k}$
is either homeomorphic to a fibration (when the map is locally
topologically  equivalent to a projection) or else it is a genuine
Montgomery-Samelson fibration with non-empty
critical locus $A$. 

If $(m,n)=(4,3)$ then the claim follows from Proposition \ref{act}. 
Let us analyze now  the case $(m,n)=(8,5)$.  
Antonelli proved in \cite{Ant3} that a 
$2k$-manifold $M$ admitting a Montgomery-Samelson 
fibration over $S^{k+1}$ is a $(k-1)$-connected $\pi$-manifold. 
It is known (see \cite{Kos2}, Proposition X.3.7, p.205)  
that  for even $k\geq 3$  any closed $(k-1)$-connected 
$\pi$-manifold $M^{2k}$ is diffeomorphic to $\Sigma^{2n}\#_{r}(S^{k}\times S^{k})$.

At last, if $M^8$ is homeomorphic to a fibration the arguments from
\cite{AF1} show that $M^8=\Sigma^{8}\# {\widehat N}$, where
$\widehat{N}$ is a fibration  over $S^{5}$ and $\Sigma^{8}$ is an
exotic sphere, as claimed.
This settles Proposition \ref{char}. \hspace{14cm}$\Box$

\begin{remark}
Kosinski  (\cite{Kos1}) also proved  
that $\Sigma^{2k}_1\#_{r}(S^{k}\times S^{k})$ 
is not diffeomorphic to $\Sigma^{2k}_2\#_{r}(S^{k}\times S^{k})$ unless 
$\Sigma_1$ is diffeomorphic to $\Sigma_2$ so that we obtain distinct 
smooth structures. 
\end{remark}

\subsection{Splitting off fibrations in dimension $(8,5)$}
Proposition \ref{char} says 
that there is essentially a unique 
way to find simply connected manifolds with finite $\varphi(M^m,S^{n})$, 
in dimensions $(4,3)$ and $(8,5)$. Manifolds $M^m$ endowed 
with some proper smooth map $M^m\to D^{n}$ 
having only finitely many critical points  
will be called {\em disk blocks}, by analogy with the spherical 
blocks above. The topology of disk blocks 
represents the essential part of the structure of pairs  with 
finite $\varphi(M^m,N^{n})$.

\begin{remark}\label{pi3}
Isomorphism classes of $S^3$-fibrations over $S^4$ are classified by 
the elements of $\pi_3({\rm Homeo}^+(S^3))$. Cerf's Theorem (\cite{Cerf}) and 
Hatcher's proof of  Smale's Conjecture (\cite{Hat}) yield  
$\pi_3({\rm Homeo}^+(S^3))\cong \pi_3(SO(4))$. It is standard 
that $SO(4)$ is diffeomorphic to $SO(3)\times S^3$ and 
$\pi_3(SO(4))\cong \pi_3(SO(3)\times S^3)\cong \Z r\oplus \Z s$, where 
the generators $r$ and $s$ are the following:
\[ r:S^3\to SO(4), r(x)y=xyx^{-1}\] 
\[ s:S^3\to SO(4), s(x)y=xy \]
and $S^3=\{ x\in {\mathbb H}, |x|=1\}$ 
is identified with set of quaternions of  unit norm. 
Here $\mathbb H$ denotes the set of quaternions. 
We consider here the same generators $(r,s)$  as those used by James and 
Whitehead in \cite{JW}, which are different from the generators 
$(h,j)$ used in \cite{Mil0}. If we denote by $\eta[a,b]$ the fibration 
associated to $ar+bs$ and by $\xi_{\alpha,\beta}$ the fibration 
associated to $\alpha h +\beta j$ then $\eta[a,b]$ is the same as 
$\xi_{\alpha,\beta}$ when 
$\alpha+\beta=b$ and $\beta=-a$. Moreover, as it is well-known, these two 
invariants are the same as the classical invariants of a rank 4 
vector bundle on $S^4$, namely the Euler class $e$  and the Pontryaguin  
class $p_1$. Specifically, we have 
\[ p_1(\eta[a,b])=(2b+4a)[S], e(\eta[a,b])=b [S]\]
where $[S]$ is the generator of $H^4(S^4)$. 
The Hopf fibration is $\eta[0,1]=\xi_{1,0}$ and has  non-zero Euler class.   
\end{remark}

Let us start now from a smooth map $f:M^m\to N^{n}$ between 
compact manifolds, having only finitely many critical 
points. We notice first that:

\begin{lemma}
The map $f$ is a generalized fiber sum $g\oplus_{\alpha} h$ of a disk 
block and a fibration over $N^n$. 
\end{lemma}
\begin{proof}
Let $D^{n}\subset N^{n}$ be a disk containing the critical values in its interior and $X^m=f^{-1}(D^{n})$. Then the restriction $f|_{M^m\setminus X^m}:M^m\setminus X^m\to N^n\setminus D^{n}$ 
is a fibration and $f$ splits as the generalized fiber sum 
of $f|_{M^m\setminus X^m}$ and $f|_{X^m}$.  
\end{proof}

Consider now that $(m,n)=(8,5)$. 
If the disk block is empty or topologically trivial it follows that the 
manifold $M^8$ is a fibration up to the connect sum with a homotopy sphere. 
Assume from now on that the disk block $h$ is non-trivial. 
Therefore, by the local structure of singularities 
(see Theorem \ref{local}) $h$ is a Montgomery-Samelson 
fibration with at least one critical point. Moreover, all its critical points 
are locally modelled by the cone over the Hopf fibration.

Recall that the double of a manifold  $X$ with boundary is the 
union of two copies of $X$ with their boundaries identified.

\begin{proposition}\label{bdhopf}
Let $h:X^8\to D^{5}$ be a smooth proper map such that 
$h|_{\partial X}$ is a fibration over $S^{4}$. We denote by 
$h\oplus h:X^8\cup_{\partial X} X^8 \to S^{5}$ the double of $h$, 
namely the map induced from the double of $X^8$ to the sphere. 
Let $X_1^8,X_2^8,\ldots, X_d^8, \ldots X_{d+f}^8$ be the connected components of $X$. 
Assume that each $X_i^8$ for $1\leq i\leq d$  contains at least one critical 
point, while those for $d+1\leq i\leq d+f$ contain none. 

Then, for each connected component $X_i^8$ with $1\leq i\leq d$  
the boundary fibration 
$f|_{\partial X_i}:\partial X_i^8\to S^4$ is a fiber sum of 
Hopf fibrations (possibly trivial) and thus is isomorphic to 
$\eta[0,r]$, for some $r\in \Z$.   
\end{proposition}
\begin{proof}
According to the local structure the map $f:X_i^8\to D^5$ has $s$
critical points $p_j$, $j=1,\ldots, s$ whose topological local model 
is the suspension of the Hopf fibration. Thus there exist local charts 
$U_j^8$ around $p_j$ and $V_j^5$ around $q_j=H(p_j)$ such that $f:U_j^8\to V_j^5$
is the cone over the Hopf fibration $h:\partial U_j\to \partial V_j$. 
Thus $U_j^8$ (and $V_j^5$) are 
8-dimensional (respectively 5-dimensional) disks. 
Moreover, the map $f: X_i^8\setminus \cup_{j=1}^s U_j^8\to D^5\setminus\cup_{j=1}^s V_j^5$ is a topological fibration.  

As in \cite{AF1}, one of the disks $U_j^8$ 
could be removed and  then glued back  in order to obtain the 
homotopy sphere summand $\Sigma^8$ and that the restriction of $f$  
to $X_i^8\setminus \cup_{j=1}^s U_j^8$ is a smooth fibration.

The fibration on each boundary component $\partial U_j$ has fiber $S^3$, and 
thus the fiber of the restriction of $f$ should be a disjoint 
union of $S^3$'s. However, the base $D^5\setminus \cup_{j=1}^s V_j^5$ is 
simply connected and thus  $X_i^8\setminus \cup_{j=1}^s U_j^8$ has as many 
connected components as the fiber. This implies that the fiber should be 
$S^3$.

\begin{lemma}\label{cobound}
Let $\alpha_1,\alpha_2,\ldots, \alpha_{s+1}$ be isomorphism classes of 
$S^3$-fibrations over $S^4$. Then there exists a $S^3$-fibration over 
$S^5\setminus\cup_{j=1}^{s+1} D^5_j$, whose restriction to the boundary 
$\partial D^5_j$ of each 5-disk is the class $\alpha_j$ if and only if 
\[ \alpha_1+\alpha_2+\cdots +\alpha_{s+1} =0 
\in \pi_3(SO(4))\cong \Z\oplus \Z\]
\end{lemma}
\begin{proof}
The holed sphere $S^5\setminus\cup_{j=1}^{s+1} D^5_j$ retracts onto 
the wedge of $s$ spheres $\vee_{\lambda=1}^{s} S^4_j$, which is 
embedded in $S^5\setminus\cup_{j=1}^{s+1} D^5_j$. Let $\alpha$ be the class 
of a $S^3$-fibration over the holed 5-sphere. We associate to $\alpha$ the 
collection $(\eta_{\lambda})_{1\leq \lambda\leq s}$ 
of restrictions to each factor $S^4_{\lambda}$ of the wedge of spheres. This 
yields an injective map from the set of isomorphism classes $\alpha$  
into  the group $\oplus_{\lambda=1}^{s} \pi_3(SO(4))$. 
This map is also surjective. Given an element $\gamma$ of 
$\oplus_{\lambda=1}^{s} \pi_3(SO(4))$ 
we consider trivial $S^3$-fibrations  on 
the 4-disks $D_{\lambda,+}^4$, $\lambda=1,\ldots, s$ 
and, respectively,  on the wedge of disks 
$\vee_{\lambda=1}^{s} D^4_{j,-}$. 
We glue together these fibrations along their boundaries, namely  
for each $\lambda$,  $\partial D_{\lambda,+}^4$ is identified to 
$\partial D_{\lambda,-}^4$ 
by using the twist $\gamma_{\lambda} \in\pi_3(SO(4))$. This yields a 
fibration over  $\vee_{\lambda=1}^{s} S^4_j$ which extends to  
$S^5\setminus\cup_{j=1}^{s+1} D^5_j$. 

We have to relate now the classes $\eta_{\lambda}$ and the classes $\alpha_j$ 
of the restrictions of $\alpha$ to the boundary components $\partial D^5_j$. 
We choose a particular embedding of  $\vee_{\lambda=1}^{s} S^4_j$ into 
$S^5\setminus\cup_{j=1}^{s+1} D^5_j$ as follows:
\begin{itemize}
\item $S_1^4$ and $\partial D_1^5$ bound a cylinder; 
\item $S_{\lambda}^4 \vee S_{\lambda+1}^4 $ and $\partial D_{\lambda+1}^5$ 
bound a 5-sphere deprived of three disks, among which two disks  
intersect  in a  boundary point, for $\lambda=1, 2, \ldots, s$; 
\item  $S_{s}^4$ and $\partial D_{s+1}^5$ bound a cylinder. 
\end{itemize}
Taking into account the orientations on each boundary component it follows 
that 
\[ \alpha_1=\eta_1,  
\alpha_j=\eta_j-\eta_{j-1},\,\; {\rm for } \,\;  2\leq j\leq s, \;\;
\alpha_{s+1}=\eta_{s}\]
Thus the classes $(\alpha_j)_{1\leq j\leq s+1}$ extend to the holed 5-sphere 
if and only if their sum vanishes. 
\end{proof}

We resume the proof of Proposition \ref{bdhopf}. We have a  fibration
$f: X_i\setminus \cup_{j=1}^s U_j\to D^5\setminus\cup_{j=1}^{s} V_j$
and all boundary fibrations over $\partial V_j$  
are Hopf fibrations. Let us denote by $\alpha_j$ the class of the 
boundary Hopf fibration $h:\partial U_j\to \partial V_j$. We   
identify each factor $\pi_3(SO(4))$ with $\Z\oplus \Z$ 
by taking account the boundary orientation of the fibrations.  
Thus each $\alpha_j$ is either the Hopf fibration, or its negative. 

For example, the suspension of the Hopf fibration $S^8\to S^5$ 
has two critical points and the classes of the associated boundary fibrations 
are the positive and the negative Hopf fibrations. 

Since the classes  
$\alpha_j$ and $f|_{\partial X_i}$ bound it follows from the 
Lemma \ref{cobound} 
that the  isomorphism class of the $S^3$-fibration $f|_{\partial X_i}$
is the (algebraic) sum of a number of Hopf fibrations and thus 
isomorphic to $\eta[0,r]$ for some integer $r$. 
It is easy to see that this is the 
same as $r$ fiber sums of Hopf fibrations. 
\end{proof}

\begin{lemma}
The connected components $X_i^8$, for $d+1\leq i\leq d+f$ are 
diffeomorphic to $D^5\times S^3$ and the fibrations 
$f|_{\partial X_i}$ are trivial $S^3$-fibrations. 
\end{lemma}
\begin{proof}
The restrictions $f:X_i^8\to D^5$ are fibrations 
whenever $d+1\leq i\leq d+f$ and we set 
$F_i^3$ for their respective fibers. 
We also set $F_i^3=S^3$ for $i\leq d$. 
According to our assumptions  the restriction of 
$f$ is also a fibration in a neighborhood of $\sqcup_{i=1}^{d+f}\partial X_i$
and thus its fiber is $\sqcup_{i=1}^{d+f}F_i$. 

Recall that the restriction to the fibrewise 
summand $f|_{M^8\setminus X^8}:M^8\setminus X^8\to N^5\setminus D^5$ is also a 
fibration and that $M^8$ is connected. The fiber on the boundary 
is known so that the fiber of $f|_{M^8\setminus X^8}$ is $\sqcup_{i=1}^{d+f}F_i^3$. 
If there exist a fiber $F_i^3$ which is not diffeomorphic to $S^3$ 
then the monodromy cannot exchange the respective connected 
components of the fiber and thus $M^8\setminus X^8$ is not connected. 
This would imply that $M^8$ is not connected, since all disk blocks 
$X_i^8$ are connected. This contradiction shows that all $F_i^3$ are diffeomorphic 
to $S^3$ and the claim follows. 
\end{proof}

\begin{lemma}\label{ext}
Let $g:M^8\setminus X^8\to N^5\setminus D^5$ be a fibration with fiber 
$\sqcup_{i=1}^{r}S^3$. Then the restrictions of $g$ to each 
boundary component is a trivial fibration and the 
fibration $g$ factors as follows: there exists a non-ramified covering 
$\widehat{N^5}\setminus \sqcup_{i=1}^r D_i^5$ of $N^5\setminus D^5$
of degree $r$ and an $S^3$-fibration $M^8\setminus X^8\to 
\widehat{N^5}\setminus \sqcup_{i=1}^r D_i^5$ and their composition is 
$g$. 
\end{lemma}
\begin{proof}
The main idea is that a $S^3$ fibration over $S^4$ extends over 
$N^5\setminus D^5$ if and only if it is trivial. This is a simple 
application of the following well-known facts: 
these invariants of the fibration can 
also be interpreted as Pontryaguin and Euler numbers of the manifold 
and the latter are cobordism invariants. 

In order to classify  (orientable) fibrations with fiber 
$\sqcup_{i=1}^{r}S^3$ one has to consider the 
group of homeomorphisms of $\sqcup_{i=1}^{r}S^3$. 
As a consequence of the Cerf and Hatcher Theorems 
this homeomorphisms group has the homotopy type of 
$SO(4)^r\ltimes S_r$, where $S_r$ denotes the 
permutation groups on $r$ elements acting on $SO(4)^r$ by permuting 
the factors. In particular isomorphism classes of such fibrations 
over $N$ are classified by the elements of the set of homotopy classes 
$[N,B(SO(4)^r\ltimes S_r)]$. 

The exact sequences between Lie groups induce the following 
exact sequence in homotopy: 
\[[N,BSO(4)]^r \to [N,B(SO(4)^r\ltimes S_r)]\to [N,BS_r]\]
Then $[N,BS_r]$ is isomorphic to the set of 
homomorphisms ${\rm Hom}(\pi_1(N), S_r)$, by associating to each 
finite covering its monodromy homomorphism. 
Moreover, given a fibration with fiber $\sqcup_{i=1}^{r}S^3$ we can associate 
a monodromy covering $\widehat{N}$ of $N$ by 
crushing each sphere to a point, or equivalently, 
by considering the monodromy action on the set of components of the fiber.

This can be interpreted as follows. 
Given a homomorphism 
$\rho:\pi_1(N)\to S_r$, the isomorphism classes 
of fibrations with fiber $\sqcup_{i=1}^{r}S^3$ 
and monodromy covering associated to $\rho$  is in bijection with 
$[N,BSO(4)]^r$. Furthermore each factor 
$[N,BSO(4)]$  contains a Euler class $e$ and a 
Pontryaguin class $p_1$, which are pull-backs of the corresponding 
elements in $K(\Z,n)$ and $K(\Z,4)$.

It is known that $\pi_j(BSO(4))=\pi_{j-1}(SO(4))$  (see e.g. \cite{St}, 
chap 19) and thus 
\[ \pi_1(BSO(4))=\pi_3(BSO(4))=0, 
\pi_2(BSO(4))=\Z/2\Z, \pi_4(BSO(4))=\Z\oplus \Z\]

By general facts concerning the Moore-Postnikov decomposition of $BSO(4)$
(see \cite{Sp}, chap. 8) we have a map 
\[ \theta: [N, BSO(4)]\to [N, K(\Z/2\Z, 2)\times K(\Z, 4)\times K(\Z, 4)]\]
or equivalently 
\[\theta:  [N, BSO(4)]\to H^2(N,\Z/2\Z)\oplus H^4(N,\Z)\oplus H^4(N,\Z)\]
When $N$ is of dimension 4 the components of $\theta(\eta)$  are  
the  characteristic classes of the fibration $\eta$, namely 
the Stiefel-Whitney class $w_2\in H^2(N, \Z/2\Z)$, the Pontryaguin class $p_1\in 
H^4(N,\Z)$ and the 
Euler class $e\in H^4(N,\Z)$. Let us denote by $e\in H^4(N,\Z)$ the corresponding 
element $\theta(\eta)$ also when $N$ is of arbitrary dimension, namely 
the image in cohomology of the pull-back of the 
generator of $H^4(K(\Z, 4),\Z)$. Notice that $e$ is not an Euler class 
if the dimension of $N$ is not $4$.

Assume now that we have a fibration $\alpha$ with fiber $\sqcup_{i=1}^{r}S^3$ 
over $N\setminus D$ and fixed monodromy homomorphism 
$\rho:\pi_1(N\setminus D)\to S_r$. We want to compute the 
invariants  of its restriction $\alpha|_{\partial D}$.
The set of isomorphisms corresponds bijectively to elements 
in $[N\setminus D, BSO(4)]^r$ and thus we have characteristic classes 
$c=(c_i)_{i=1,r}\in \oplus_1^r H^4(M)$ corresponding to each one 
of the $r$ factors.

For any characteristic class $c\in H^4(M)$ we have  
\[ c(\alpha|_{\partial D})=i^*c(\eta)\]
where $c$, $i:\partial D\to N\setminus D$ 
is the inclusion and $i^*:H^4(N\setminus D)\to H^4(\partial D)$ 
is the map induced in cohomology. 
However the map $i^*$ is trivial. For instance   in the Mayer-Vietoris 
sequence 
\[ 0\to H^4(N)\to H^4(N\setminus D)\oplus H^4(D^5)\to H^4(\partial D)\to \] 
the map $H^4(N)\to H^4(N\setminus D)$ is an isomorphism and thus 
the kernel of $i^*$ is all of $ H^4(N\setminus D)$. 
Thus the  characteristic classes of 
each $S^3$-fibration component of 
$\alpha|_{\partial D}$ vanish. 

This  result holds more generally for any boundary $4$-manifold 
instead of the sphere $\partial D$, 
since according to Dold and Whitney (see \cite{DW, Woo}) 
the classes $w_2,p_1, e$ classify $S^3$-fibrations over 4-manifolds. 

By Remark \ref{pi3}  each component of 
$\alpha|_{\partial D}$ is a trivial fibration. 
\end{proof}

\begin{lemma}
Each connected component $X_i^8$ with $1\leq i\leq d$ is diffeomorphic to 
some manifold $\Sigma_i^8\#_{r_i}S^4\times S^4 \setminus D^5\times S^3$, 
so that the restriction of $f$ to the boundary is  identified with the projection 
$S^3\times \partial D^5\to \partial D^5$ on the second factor. 
\end{lemma}
\begin{proof}
The boundary fibration $f|_{\partial X_i}$ is a trivial fibration, and hence 
we can glue to $X_i^8$ the trivial fibration over $D^5$ in order to obtain 
a manifold $Y_i^8$ endowed with a smooth map with only finitely many 
critical points into $S^5$. 

Further $\partial X_i$  and $X_i\setminus \cup_{j=1}^s U_j$, 
and hence $X_i$, are simply connected as $S^3$-fibrations 
over simply connected bases. Thus the manifold 
$Y_i$ is a  simply connected  spherical block. 
The structure of  spherical blocks from 
Proposition \ref{char} implies that $Y_i^8$ is diffeomorphic to 
$\Sigma_i^{8}\#_{r_i} S^{4}\times S^{4}$ for some $r_i\geq 0$, with 
$\Sigma_i^8$ a
homotopy sphere. 
\end{proof}





\begin{remark}
The  
isomorphism classes of $S^3$-fibrations over  
$\#_c S^1\times S^4$ are classified by 
an element in $\oplus_1^c\pi_3(SO(4))=(\Z\oplus\Z)^c$, 
each factor $\Z\oplus \Z$ corresponding to the isomorphism class of the 
restriction of the $S^3$-fibration to a factor $\{*\}\times S^4$. 
\end{remark}

\subsection{Splitting off fibrations in dimensions $(4,3)$}
The aim of this section is to prove the corresponding splitting result 
in dimension $(4,3)$. Some extra care is needed to 
describe the smooth structure of the 4-manifolds. 

The first step in proving the Theorem is to understand the smooth 
structure of spherical blocks in dimension $(4,3)$, by a slight strengthening 
of Proposition \ref{char}: 

\begin{proposition}\label{sblock4}
If $M^4\to S^3$ has  finitely many critical points 
and its generic fiber is $S^1$, then 
\begin{enumerate}
\item either $M^4=\Sigma^4\# N^4$, where $N^4$ is a fibration over $S^4$;
\item or else $M^4$ is $\Sigma^4\#_r S^2\times S^2$, 
where 
$\Sigma^4$ is a homotopy 4-sphere. 
\end{enumerate}
\end{proposition}
\begin{remark}
If $\pi_1(M^4)=0$ then the generic fiber should be a circle.  
It will be proved later (see Proposition \ref{subgroup}) that a connected 
spherical block  $M^4$ which is non-fibered is simply connected. 
On the other hand a fibered  connected spherical block $M^4$ 
is a circle fibration over $S^3$. 
Thus, the generic fiber is a circle for any 
connected spherical block.     
\end{remark}
\begin{proof}
By the proof of Proposition \ref{char} 
a smooth map $f:M^4\to S^3$ with finitely many 
critical points is a Montgomery-Samelson fibration. Thus $f$ is locally topologically equivalent to either the standard projection or else to the 
suspension $H$ of the Hopf fibration $S^3\to S^2$ around branch points.  
We will use a refinement of the argument from \cite{AF1} 
for converting a  topological fibration into a smooth fibration at the 
expense of adjoining some homotopy sphere.  

Around each critical point $p\in M^4$ from the branch locus there 
exists an open neighborhood $p\in U^4\subset M^4$ 
and a homeomorphism $g:V^4 \to U^4$ such that 
$f\circ g:V^4\to f(U^4)$ is smoothly equivalent to the suspension $H$. 
Moreover, by choosing a smaller neighborhood $V^4$, we can assume 
that $V^4$ is a 4-disk and that its image $U^4$ is contained in a 4-disk 
inside $M^4$. This implies that $U^4$ is a homotopy 4-disk. 
If we remove the homotopy 4-disk $U^4$ and glue back 
$V^4$ instead then we obtain the manifold $M^4\# \Sigma_1^4$, where 
$\Sigma_1^4$ is a homotopy sphere, namely the negative (in the group 
of homotopy 4-spheres) of the homotopy sphere 
obtained by capping off the boundary of $U^4$ by a 4-disk.
A Theorem of Huebsch and Morse  (see \cite{HM}) says that 
whenever $\Delta^4$ is a homotopy 4-disk which can be embedded 
in a 4-disk there exists a smooth 
homeomorphism $\Delta^4\to D^4$ with only one critical point.  
As a consequence there exists  a smooth 
homeomorphism $U^4\to V^4$ having precisely one critical point.

It follows that there is an induced smooth map  $M^4\# \Sigma_1^4\to S^3$, 
which coincides with $f$ outside $V^4$ and it  is locally smoothly 
equivalent around the image of $p$ to the suspension map $H$. 
We continue in the same way for all points of the branch locus and 
get a manifold $M\# \Sigma^4$ endowed with a map to $S^3$ 
whose branch locus critical points are locally smoothly equivalent to 
 $H$ and $\Sigma^4$ is the connected sum of the homotopy 4-spheres 
corresponding to each critical point. 

At last, consider the critical points of $f$ which are not in the 
branch locus. Locally the map $f$ is topologically equivalent to the 
projection. As above, by adjoining a homotopy 4-sphere we can make the 
map be smoothly equivalent to a projection and hence a submersion. 

At the end of this process we obtained a smooth map 
$f:N^4=\Sigma^4\# M^4\to S^3$ locally modelled on $H$. 
If the branch locus is empty then $f$ is a fibration. 

Assume that the set of branch points is non-empty from now on.

Further $S^1$-fibrations  over $S^2$ are classified by 
an element of $\pi_1({\rm Homeo^+}(S^1))\cong \pi_1(SO(2))=
\Z$, namely an integer 
called the Euler class of the fibration. 

Lemma \ref{cobound} has a similar statement in this dimension, 
as follows: 

\begin{lemma}\label{cocobound}
Let $\alpha_1,\alpha_2,\ldots, \alpha_N$ be isomorphism classes of 
$S^1$-fibrations over $S^2$. Then there exists a $S^1$-fibration over 
$S^3\setminus\cup_{j=1}^N D^3_j$, whose restriction to the boundary 
$\partial D^3_j$ of each 5-disk is the class $\alpha_j$ if and only if 
\[ \alpha_1+\alpha_2+\cdots +\alpha_N =0 \in \pi_1(SO(2))\cong \Z\]
\end{lemma}
\begin{proof}
We omit it. 
\end{proof}

If we excise disk neighborhoods around the $s$ critical values of the map 
$f:N^4\to S^3$ we obtain a $S^1$-fibration over 
$S^3\setminus\cup_{j=1}^{s} D^3_j$ extending the boundary (signed) 
Hopf fibrations. Since the Euler class of the 
Hopf fibration is the generator $1\in \pi_1(SO(2))$, and thus  it is non-zero, 
Lemma \ref{cocobound} implies that $s=2k$ is even and there 
exist $k$ positive and $k$ negative Hopf fibrations.

Let consider embedded 2-spheres $S_1^2,\ldots,S_{k-1}^2$ separating 
boundary components into pairs of opposite signs.  
Using again Lemma \ref{cocobound} a number of times we obtain that 
the restriction of the $S^1$-fibration $f$ at $S_j^2$ is a trivial fibration, 
for each $j$. This means that the fibration $f$ is the fiber sum 
of spherical blocks with 2 critical values along trivial fibrations. 

\begin{lemma}
If $f:N^4\to S^3$ is as above and  has two critical values then $N^4$ is 
diffeomorphic to $S^4$.   
\end{lemma}
\begin{proof}
We have  disk neighborhoods $D^4_+$ and $D^4_-$ in $N^4$
around the critical points and $f:N^4\setminus (D^4_+\cup  D^4_-)\to S^2\times [0,1]$ 
is a $S^1$-fibration. The restrictions of $f$ to the 
boundary spheres $\partial D^4_+$ and $\partial D^4_-$ are isomorphic to the 
Hopf fibrations over $S^2$.  The composition of the projection   
$S^2\times [0,1]\to [0,1]$ with the above restriction of $f$  is still a fibration of 
$N^4\setminus (D^4_+\cup  D^4_-)$ onto $[0,1]$, and hence 
a trivial one.  Therefore  $N^4\setminus (D^4_+\cup  D^4_-)$ is diffeomorphic 
to $S^{3}\times [0,1]$. 
The manifold $N^4$ is obtained by gluing $D^4_+$, $S^{3}\times [0,1]$ and 
$D^4_-$ along their respective boundaries. Cerf's Theorem $\Gamma_4=0$ 
(see \cite{Cerf}) implies that $N^4$ is diffeomorphic to $S^4$, as claimed. 
\end{proof}

Let now $M^4_r$ denote the total space of the fiber sum of $r$ blocks
$S^4\to S^3$. Thus $M^4_1=S^4$ and $M^4_{r+1}=M^4_r\oplus S^4$, for $r\geq 1$. 
We obtain $M^4_{r+1}$ by deleting out neighborhoods $S^1\times {\rm int}(D^3)$ 
of generic fibers from $M^4_r$ and $S^4$ and gluing then 
$M^4_r\setminus S^1\times {\rm int}(D^3)$ and $S^4\setminus S^1\times {\rm int}(D^3)$  
along their boundaries, by a gluing homeomorphism  respecting 
the product structure on the boundary.  
 
Consider  first the case $r=2$. Since the fiber $S^1\subset S^4$ is unknotted it follows that 
$S^4\setminus S^1\times {\rm int}(D^3)$ is diffeomorphic to 
$S^2\times D^2$. Moreover we want  to glue the two copies of $S^2\times D^2$ 
such that the two boundary fibrations  $S^2\times S^1\to S^2$ 
glue together, i.e. 
the gluing homeomorphisms $\Phi:S^2\times S^1\to S^2\times S^1$ induces a homeomorphism 
of the base $\varphi:S^2\to S^2$ so that the diagram 
\[ \begin{array}{ccc}
S^2\times S^1 & \stackrel{\Phi}{\to}    &  S^2\times S^1\\
\downarrow     &                                     & \downarrow \\
S^2                 & \stackrel{\varphi}{\to} & S^2\\
\end{array}
\]
is commutative. 
This condition implies that the result of the gluing of the two spherical blocks 
possesses a smooth map into $D^3\cup_{\varphi} D^3$, which is a sphere  $S^3$, since 
diffeomorphisms reversing the orientation of $S^2$ are isotopic to  
a symmetry, by a classical result of Smale. 

The isotopy class of a gluing homeomorphism  respecting the product structure corresponds to 
a homotopy class of a map $S^2\to {\rm Homeo}^+(S^1)$, 
namely  to an element of 
$\pi_2(SO(2))=0$. 
Therefore the diffeomorphism type of $M_2$  does not depend 
on the choice of this isomorphism.

Therefore $M_2$ is the result of gluing  two copies of 
$S^2\times D^2$ by using the boundary identification $\Phi$. 
Observe that there is also an  obvious projection 
$S^2\times D^2\to S^2$ extending the boundary projection
$S^2\times S^1\to S^2$. The above diagram shows that the 
two projections on the first factor of each copy $S^2\times D^2$  
glue together to get a trivial fibration map $M^4_2\to S^2$. 
Moreover,  the fiber of this fibration is the 
union $D^2\cup D^2$, namely a 2-sphere $S^2$. 
It follows then that $M_2=S^2\times S^2$.

When $r>2$, we should also analyze what happens from $S^4$ when we delete 
neighborhoods of two fibers. Since the fibers form a trivial link in $S^4$  and can be separated 
by an embedded 3-sphere it follows that 
$S^4\setminus (S^1\times {\rm int}(D^3)\cup S^1\times {\rm int}(D^3))$ is diffeomorphic 
to the connected sum  out of the boundary 
$(S^4\setminus S^1\times {\rm int}(D^3))\# 
(S^4\setminus S^1\times {\rm int}(D^3))$. 
Recall that $M_r$ is obtained by  iterating boundary gluings of 
$S^4\setminus S^1\times {\rm int}(D^3)$  with $(r-2)$ copies of  
$S^4\setminus (S^1\times {\rm int}(D^3)\cup S^1\times {\rm int}(D^3))$ 
glued each one to the next one and finally with one more copy of 
$S^4\setminus S^1\times {\rm int}(D^3)$. 
By the previous observation this gluing is the connected sum of $(r-1)$ manifolds, where 
each manifold is  obtained by gluing two $S^4\setminus S^1\times {\rm int}(D^3)$. 
It follows then than $M_r=\#_{r-1}S^2\times S^2$. 
This proves Proposition \ref{sblock4}. 
\end{proof}

\subsection{End of proof of Theorem \ref{complete}}

\subsubsection{Dimension (8,5)}
The fibrewise summand $f|_{M^8\setminus X^8}$ extends a trivial 
boundary fibration. We can glue then along the boundary a trivial fibration 
over the disk $D^5$ to obtain a fibration $g:W^8\to N^5$ with fiber
$\sqcup_1^{d+f}S^3$. This implies that $W^8$ is a $S^3$-fibration over the 
monodromy covering $\widehat{N}$, which is a degree $(d+f)$ covering 
of $N$. 

Since the disk blocks over $D_i^5$ are connected and disjoint, the gluing of  
$M^8\setminus X^8$ and $X^8$  to recover $M^8$ corresponds to taking 
$d+f$ fiber sums (along 
$d+f$  components $S^3$, which are fibers of $W^8$) with the manifolds $Y_i^8$. 
We put also $Y_j^8=S^3\times S^5$, when $d+1\leq j\leq d+f$. 

Consider now $\Sigma^8=\Sigma_1^8\# \cdots \# \Sigma_d^8$, where $\Sigma_j^8$ are 
the homotopy 8-spheres associated to each factor. 
Then we can slide the attaching maps of  $\Sigma_i^8$  in order to be attached  
to a small disk on $W^8$. This is possible since the image of the 
monodromy  of the covering is all of the permutation group 
$S_{d+f}$, as $M^8\setminus X^8$ is connected. 
We obtained the structure claimed in the Theorem for 
$D=d+f$, where all factors $Y_i^8$ (with $i\geq d+1$) are diffeomorphic to 
$\#_{-1}S^2\times S^2$.

Let $q=2\sum_{i=1}^dr_i+2d$. 
Note that $\varphi(M^8, N^5)\leq q$, because 
we have already a smooth map with $q$ critical points.

\subsubsection{Dimension (4,3)}
Consider the disk block $X^4\to D^3$ and assume that $X^4$ has $d+f$ connected 
components $X_i^4$ such that the first $d$ have at least one branch point. 
Then $X_i^4$ minus a number of spherical caps is a fibration 
over the holed $D^3$. Since the holed 3-disk is simply connected it follows 
that the  generic fiber of $f|_{X_i^4}$ is a single $S^1$.

Therefore $M^4\setminus X^4$ is a manifold with $d+f$ boundary 
components, which are fibrations over $\partial D$.

\begin{lemma}\label{ext4}
Let $g:M^4\setminus X^4\to N^3\setminus D^3$ be a fibration with fiber 
$\sqcup_{i=1}^{r}S^1$. Then the restrictions of $g$ to each 
boundary component is a trivial fibration and the 
fibration $g$ factors as follows: there exists a non-ramified covering 
$\widehat{N^3}\setminus \sqcup_{i=1}^r D_i^3$ of $N^3\setminus D^3$
of degree $r$ and an $S^1$-fibration $M^4\setminus X^4\to 
\widehat{N^3}\setminus \sqcup_{i=1}^r D_i^3$ and their composition is 
$g$. 
\end{lemma}
\begin{proof}
We can follow the same lines as in the proof of Lemma \ref{ext} by replacing 
$BSO(4)$ with $BSO(2)$. However there is a much simpler proof. 
In fact a circle fibration over $\partial D^3$ which extends 
to $N^3\setminus D^3$ is trivial because 3-manifolds are parallelizable.
We leave the details to the reader.  
\end{proof}

\begin{lemma}
Each connected component $X_i^4$ with $1\leq i\leq d$ is diffeomorphic to 
some manifold $\Sigma_i^4\#_{r_i}S^2\times S^2 \setminus D^3\times S^1$,
so that the restriction of $f$ to the boundary is the projection 
$S^1\times \partial D^3\to \partial D^3$. 
\end{lemma}
\begin{proof}
We construct the manifold $Y_i^4$ by adjoining to $X_i^4$ 
a trivial $S^1$-fibration over 
$D^3$. Although we cannot insure that $Y_i^4$ is simply connected we know that 
$Y_i^4$ admits a smooth map with finitely many critical points into 
$S^3$ having the generic fiber $S^1$. Thus we can apply Proposition 
\ref{sblock4} to find that $Y_i^4$ is diffeomorphic to $\Sigma_i^4\#_{r_i}S^2\times S^2$. 
\end{proof}

The proof of the Theorem follows from these Lemmas.

\section{Computing $\varphi$ and the proof of Theorem \ref{compute}}

Computing the precise value for $\varphi$ from the structure 
Theorem \ref{complete} is easy when the fundamental group is trivial 
or verifies some strong assumptions forcing the  finite 
covering to be unique. In fact the main point is to understand 
whether a given $M^m$ admits  several fiber sums decompositions 
associated to coverings $\widehat{N^n}$ of different degrees.

Henceforth we consider  two manifolds  $M^{2k},N^{k+1}$ 
with  finite $\varphi(M^{2k},N^{k+1})$ and $k\in \{2,4\}$. 
According to Theorem \ref{complete} 
either $M^{2k}$ is a topological fibration over $N^{k+1}$ 
or else  it is the iterated fiber 
sum of a fibration  $W^{2k}$  over $N^{k+1}$ 
with  a number of spherical blocks. 
In the second case we say, by a language abuse, that 
$M^{2k}$ is {\em non-fibered} (over $N^{k+1}$). 

 We have the following result (restatement of Corollary \ref{gp} 
from Introduction):   
 
\begin{proposition}\label{subgroup}
Let $M^{2k},N^{k+1}$ be closed orientable manifolds with finite 
$\varphi(M^{2k},N^{k+1})$, $k\in\{2,4\}$ and $M^{2k}$ non-fibered.  

If  $k=4$  then $\pi_1(M^8)$ is a finite index 
normal subgroup of $\pi_1(N^5)$. 

If $k=2$  then  $\pi_1(M^4)\cong\pi_1(N^3)$.    
\end{proposition}
\begin{proof}
Let $F^{k-1}$ denote the fiber of $W^{2k}\to N^{k+1}$.
 
If $k=4$ then $\pi_1(W^8)\cong \pi_1(M^8)$ by Van Kampen. The 
homotopy exact sequence 
of the fibration $W^8\to N^5$ yields the exact sequence: 
\[ 0\to \pi_1(M^8)\to \pi_1(N^5)\to \pi_0(F^3)\to 0\]
This means that $\pi_1(M^8)$ is a normal finite index subgroup of 
$\pi_1(N^5)$. 

If $k=2$ the homotopy exact sequence of the fibration reads 
\[ 0\to \pi_2(M^4)\to \pi_2(N^3)\to \pi_1(F^1)\to \pi_1(W^4)\to \pi_1(N^3)\to \pi_0(F^1)\to 0\]
We obtain $M^4$ from $W^4$ by capping off some of the 
$S^1$ components of the fiber. In fact, for realizing 
the fiber sum we remove a neighborhood of the fiber component and replace 
it with the complement of a neighborhood $Z^4$ of the generic fiber 
in some $\#_r S^2\times S^2$. However, a fiber in the 
boundary of $Z^4$ is null-homotopic in 
$\#_r S^2\times S^2\setminus Z^4$: translate the fiber until it lies 
in the boundary of the cone of a Hopf fibration around a critical point 
and then collapse it onto the vertex of the cone.

Since the monodromy homomorphism of the covering $\widehat{N^3}\to N^3$ at the 
level of connected components  of the fiber 
is surjective capping off one component 
is equivalent to capping off any other component at the level of 
fundamental group. Specifically, by Van Kampen 
\[ \pi_1(M^4)=\pi_1(W^4)/\langle i(\pi_1(F^1))\rangle\]
where $\langle i(\pi_1(F^1))\rangle$ is the normal subgroup of 
$\pi_1(N^3)$ generated by the image $i(\pi_1(F^1))$ by the inclusion map 
$i:F^1\hookrightarrow W^4$. 

The exact sequence above shows that the image $i(\pi_1(F^1))$ is a normal 
subgroup, namely the kernel  of $\pi_1(W^4)\to \pi_1(N^3)$. Therefore 
there exists a natural identification of $\pi_1(M^4)$ with $\pi_1(N^3)$. 
\end{proof}

\begin{definition}
The group $G$ is {\em co-Hopfian} if any injective homomorphism 
$G\to G$ is an automorphism. 

\end{definition}

\begin{proposition}\label{hopfian}
Let $M^{2k},N^{k+1}$ be closed orientable manifolds with finite 
$\varphi(M^{2k},N^{k+1})$, $k\in\{2,4\}$ and $M^{2k}$ non-fibered.  
 
If $k=4$, assume that  $\pi_1(M^8)\cong \pi_1(N^5)$ is a co-Hopfian group 
or a finitely generated free non-abelian group.  Then 
$\varphi(M^8,N^5)=2r+2D$, where $r=r_1+r_2+\cdots+r_D$.

If $k=2$,  assume that $N^3$ is an irreducible 
closed orientable 3-manifold which is not finitely covered by 
a torus bundle over $S^1$ nor by  $\Sigma\times S^1$ 
(for some surface $\Sigma)$.  
Then $\varphi(M^4,N^3)=2r+2D$. 
\end{proposition}
\begin{proof}
Let $k=4$. The structure Theorem and Van Kampen 
imply that $\pi_1(W^{8})\cong \pi_1(M^{4})$.  
Moreover, the homotopy exact sequence of the 
fibration $W^8\to N^5$ implies that the fiber is a single sphere, because 
the injective homomorphism $\pi_1(M^8)\to \pi_1(N^5)$ should be an automorphism 
if $\pi_1(N^5)$ is co-Hopfian. 

If $\pi_1(N^5)$ is free non-abelian then an argument also used in \cite{FPZ} 
can be applied. The image of $\pi_1(M^8)$ into $\pi_1(N^5)$ should be  a normal 
subgroup of finite index equal to the number of components of the fiber. 
By the Nielsen-Schreier Theorem this is a free non-abelian group of 
rank greater than the rank of $\pi_1(M^8)$ and thus cannot be the image 
of $\pi_1(M^8)$ unless the fiber is connected. 

This implies that there exists only one possible value for $D$, namely  
$D=1$. 

\begin{lemma}\label{unique}
If $D=1$  and $k\in\{2,4\}$ then $r$ is uniquely determined by 
the homology of $M^{2k}$ and 
$N^{k+1}$. Specifically, we have:  
\[ 2r = \left\{\begin{array}{ll}
b_k(M^{2k})-2b_1(N^{k+1}), & \mbox{ \rm if } b_k(M)\equiv 0({\rm mod}\, 2)\\
b_k(M^{2k})-2b_1(N^{k+1})+1, & \mbox{ \rm if } b_k(M)\equiv 1({\rm mod}\, 2)
\end{array}
\right.
\]
where $b_j$ states for the  $j$-th Betti number. 
\end{lemma}
\begin{proof}
It suffices to consider $r\geq 0$ since otherwise $M^{2k}=W^{2k}$ 
and $M^{2k}$ fibers over $N^{k+1}$. 

All homology groups below will be considered with rational 
coefficients. The Gysin sequence of the fibration 
$W^{2k}-S^{k-1}\times D^{k+1}\to N^{k+1}-D^{k+1}$  yields 
\[ 0\to H_1(N^{k+1}-D^{k+1})\to H_k(W^{2k}-S^{k-1}\times D^{k+1})\to 
H_k(N^{k+1}-D^{k+1})\stackrel{\beta}{\to} H_0(N^{k+1}-D^{k+1})\]
and thus 
\[ {\rm rk}_{\Q} H_k(W^{2k}-S^{k-1}\times D^{k+1})=
2b_1(N^{k+1})-{\rm rk}_{\Q} {\rm Im}(\beta)\]
where ${\rm rk}_{\Q}$ denotes the rank over $\Q$. 
Moreover the image of $\beta$ is contained in $H_0(N^3-D^3)$ and thus 
$0\leq {\rm rk}_{\Q} {\rm Im}(\beta)\leq 1$. 
Further $H_k(N^{k+1}-D^{k+1})\cong H_k(N^{k+1})\oplus H_k(\partial D^{k+1})$ and 
the map $\beta$ on $H_k(N^{k+1})$  is the cap product with the Euler 
class $e_{W}\in H^k(N^{k+1})$ of the fibration $W^{2k}\to N^{k+1}$ in $(k-1)$-spheres and
 respectively trivial on the factor $H_k(\partial D^{k+1})$. It follows 
that 
\[ {\rm rk}_{\Q} {\rm Im}(\beta)=\left\{\begin{array}{ll}
0, & \mbox{ if } e_{W}=0 \\
1, & \mbox{ otherwise } 
\end{array}
\right.
\]

We want now to compute the homology of $M^{2k}$ which is obtained by gluing 
$W^{2k}-S^{k-1}\times D^{k+1}$ with $Z^{2k}=\#_{r}S^k\times S^k \setminus S^{k-1}\times D^{k+1}$ 
along $S^{k-1}\times S^k$. The proof of Proposition \ref{sblock4} shows us that 
$Z^{2k}$ is diffeomorphic to 
$\#_{r}S^k\times S^k \# S^k\times D^k$ and 
thus $H_k(Z^{2k})=\Q^{2r+1}$ and $H_1(Z^{2k})=0$, if $r\geq 0$. 
By Mayer-Vietoris we derive the exact sequence:  
\[ \Q=H_2(S^1\times S^2) \stackrel{\zeta}{\to} 
H_k(W^{2k}-S^{k-1}\times D^{k+1})\oplus H_k(Z^{2k}) \to H_k(M^{2k}) \stackrel{\nu}{\to}
H_{k-1}(S^{k-1}\times S^k) =\Q\]
so that 
\[ {\rm rk}_{\Q} H_k(M^{2k})={\rm rk}_{\Q} H_k(W^{2k}-S^{k-1}\times D^{k+1}) + 2r+1 +
{\rm rk}_{\Q} {\rm Im}(\nu) -{\rm rk}_{\Q} {\rm Im}(\zeta)\]
The map $\zeta$ is injective since $H_k(Z^{2k})$ is generated by the  
homology classes of  the obvious embedded 2-spheres. 
Thus $\ker \zeta=0$ and hence ${\rm rk}_{\Q} {\rm Im}(\zeta)=1$. 

Note that 
\[ 0\leq {\rm rk}_{\Q} {\rm Im}(\nu) \leq 1\]
Observe now that the map induced by inclusion 
$H_{k-1}(S^{k-1}\times S^k)\to H_{k-1}(Z^{2k})$ should be trivial since 
the fiber $S^{k-1}\times \{*\}$ is null-homologous in $Z^{2k}$. 

We have $H_{k-1}(W^{2k}-S^{k-1}\times D^{k+1})\cong H_{k-1}(W^{2k})$ 
because any $(k-1)$-cycle can be 
homotoped  in $W^{2k}$ outside the fiber $S^{k-1}$, by general position arguments. 
The map induced by inclusion $H_{k-1}(S^{k-1}\times S^k)\to  
H_{k-1}(W^{2k}-S^{k-1}\times D^{k+1})
\cong H_{k-1}(W^{2k})$  is trivial if  and only if 
the fiber of $W^{2k}\to N^{k+1}$ is null-homologous in $W^{2k}$. 
Therefore, by Mayer-Vietoris we have  ${\rm rk}_{\Q} {\rm Im}(\nu) = 1$ if and only if the 
fiber of $W^{2k}\to N^{k+1}$ is null-homologous, and thus ${\rm rk}_{\Q} {\rm Im}(\nu) = 0$ 
otherwise. 

Now, if $e_{W}=0$ then the fiber cannot be null-homologous since the fibration  $W$ 
is trivial and thus
\[ {\rm rk}_{\Q} H_k(M^{2k})=2b_1(N^{k+1})+2r\]
When $e_{W}\neq 0$ we have 
\[ {\rm rk}_{\Q} H_k(M^{2k})=2b_1(N^{k+1})+2r-1 + {\rm rk}_{\Q} {\rm Im}(\nu) \]
In particular we need that  
${\rm rk}_{\Q} {\rm Im}(\nu)\equiv  b_k(M^{2k})({\rm mod }\, 2)$ and the 
claimed formula follows. 
\end{proof} 

\begin{remark}
Similar computations can be done for arbitrary $D$ but instead of $b_1(N)$ 
we have to use $b_1(\widehat{N})$.  However, we have 
no control on the  first Betti number of finite coverings in dimension three, 
out of the trivial upper bound. 
Thus the above result does not extend when $D>1$.   
\end{remark}

Therefore $D$ and $r$ are uniquely determined by the topology of $M^{2k}$ 
and $N^{k+1}$. 
Consider a  smooth function $M^{2k}\to N^{k+1}$ with 
$\varphi(M^{2k},N^{k+1})$ critical points. The proof of  
Theorem \ref{complete} shows that there exists some decomposition 
of $M^{2k}$ as a fiber sum $W^{2k}\oplus\#_r S^k\times S^k$ 
such that $\varphi(M^{2k},N^{k+1})=2r+2$ 
(recall that $D=1$), under the assumptions of 
Proposition \ref{hopfian}. Then Lemma \ref{unique} shows that any 
decomposition of $M^{2k}$ 
as a fiber sum $W^{2k}\oplus\#_s S^k\times S^k$, where 
$W^{2k} \to N^{k+1}$ factors as a spherical fibration 
over some degree $D$ covering of $N^{k+1}$,  actually forces $D=1$ and $s=r$. 
Thus $\varphi$ is as claimed and this ends the proof of the case $k=4$.

Let us consider now the case $k=2$. 
We can assume that the closed orientable  3-manifold 
$N^3$ is geometric, according to Perelman's solution 
to Thurston's geometrization Conjecture. 
Since $\pi_1(N^3)$ is co-Hopfian $N^3$ is irreducible. 
Moreover, if $N^3$ is not a Seifert fibered manifold covered 
by a torus bundle over $S^1$ nor by a $\Sigma\times S^1$, then 
$\pi_1(N)$ is co-Hopfian (see \cite{WW,YW}).

Let us recall that $f:M^4\to N^3$ is a Montgomery-Samelson fibration 
and thus it induces a fibration 
$M^4\setminus V\to N^3\setminus B$, where $B$ and $V$ are  the finite sets 
of critical values and  critical points, respectively. 
Then the 
homotopy exact sequence of this fibration yields 
\[ 0\to \pi_2(M\setminus V)\to \pi_2(N\setminus B)\to \pi_1(F)\to 
\pi_1(M\setminus V)\to \pi_1(N\setminus B)\to \pi_0(F)\to 0\]
Note that $\pi_1(M)\cong \pi_1(M\setminus V)$ and 
  $\pi_1(N)\cong \pi_1(N\setminus B)$. 
We observed in the proof of the Proposition \ref{subgroup} above that 
the inclusion $F\hookrightarrow M$ sends $\pi_1(F)$ to $0$.
Therefore, we derive  the exact sequence:
\[ 0\to \pi_1(M)\to \pi_1(N)\to \pi_0(F)\to 0\]
But $\pi_1(M)$ and $\pi_1(N)$ are isomorphic. Thus,  
if $\pi_1(N)$ is co-Hopfian then $\pi_0(F)$ is 0. 
This implies that the fibration  $W$ has connected fibers. 
Therefore $D=1$ and hence Lemma \ref{unique}  shows that 
$r$ is uniquely determined.
The same argument which was used above for $k=4$  ends the proof.    
\end{proof}

\begin{remark}
\begin{enumerate}
\item The  braid group $B_3$, which is the fundamental 
group of the trefoil knot,  
is a well-known non-co-Hopfian group (see \cite{Ro}). 
For instance there exist finite index subgroups of any index 
$k\equiv\pm 1({\rm mod} \,6)$ which are isomorphic to $B_3$. 
Its quotient $B_3/\Z=PSL(2,\Z)$ is also non-co-Hopfian. 
\item 
There exist also non-co-Hopfian 3-manifold groups, for instance the 
following groups: 
\[ \langle a,b,t| ab=ba, \, tat^{-1}=a^ub^v, \, tbt^{-1}=a^wb^z\rangle \]
where 
$\left(\begin{array}{cc}
u & v \\
w & z \\
\end{array}\right)$ 
is a (hyperbolic) matrix with two real distinct eigenvalues.  
These are fundamental groups of torus bundles over the circle 
which are sol-groups i.e. lattices in the group SOL. 
\end{enumerate}
\end{remark}

\begin{proposition}\label{noncohopf}
Let $M^4,N^3$ be closed orientable manifolds with finite $\varphi(M^4,N^3)$ 
and $M^4$ non-fibered. Suppose that  the fundamental group of the 
closed orientable $3$-manifold $N^3$ is not co-Hopfian.  
Then, for any decomposition of $M^4$ as a fiber sum, the value of  
$2r+2D$ is independent on the choice of the 
covering and 
\[ \varphi(M^4,N^3)=2r+2D= 2\left[\frac{b_2(M)+1}{2}\right]-2b_1(N)+2\]
\end{proposition}
\begin{proof}
Consider first $N^3$ be a closed irreducible 3-manifold. 
Assume that we have  a  finite non-trivial  covering $\widehat{N^3}$ of $N^3$ 
with  $\pi_1(\widehat{N^3})\cong\pi_1(N^3)$.  Then 
$\pi_1(N^3)$ is infinite and not co-Hopfian. 

\begin{lemma}
The covering $\widehat{N^3}$ is homeomorphic to $N^3$. 
\end{lemma}
\begin{proof}
We will use below the fact that all 3-manifolds are geometric, as it was established by 
Perelman.

The  geometric 3-manifolds whose groups  are not co-Hopfian are either 
reducible or else are finitely covered by either a 
torus bundle or else by a product $\Sigma\times S^1$, where 
$\Sigma$ is a  closed surface (see \cite{WW,YW}). 

If $N^3$ is covered  by $\Sigma\times S^1$ then it is a Seifert fibered 3-manifold. 
In particular  $N^3$ is a closed aspherical 3-manifold whose homotopy type is 
completely described by $\pi_1(N^3)$. The finite covering $\widehat{N^3}$ 
is also an irreducible Seifert fibered 3-manifold with the same fundamental group.  
The asphericity implies that $N^3$ and $\widehat{N^3}$ are homotopy equivalent. 
A Theorem of Scott (\cite{Sc}) implies that they are actually homeomorphic.

If $N^3$ is finitely covered by a torus bundle over $S^1$ then the image of the 
fiber is an incompressible surface in $N^3$, corresponding to the embedding 
of $\Z\oplus\Z$ in $\pi_1(N^3)$. Therefore $N^3$ is Haken. 
Since $N^3$ and $\widehat{N^3}$ are aspherical they are homotopy equivalent and 
a classical Theorem of Waldhausen shows that they are homeomorphic.  
\end{proof}

In particular, the fibration $W^3$ is diffeomorphic to a $S^1$-fibration over 
$N^3$. Therefore $M^3$ is diffeomorphic to the fiber sum along $D$ 
distinct fibers of the $S^1$-fibration $W^4$ with $D$  respective 
spherical blocks $\#_{r_i}S^2\times S^2$. 
Then there exists a smooth function $M^4\to N^3$  
with $2r+2D$ critical points, which has connected fibers.
Take a disk in $N^3$ containing the critical values and its 
inverse image under the smooth map, which is a connected  disk block. 
The proof of Theorem \ref{complete} shows that $M^4$ decomposes as the 
fiber sum of a spherical block containing the $2r+2D$ critical points  
and a circle fibration $V^4\to N^3$. But the spherical block is 
diffeomorphic to $\#_{(\sum_{i=1}^Dr_i)+D-1}S^2\times S^2$. 
In the new  fiber sum decomposition of $M^4$   we have  then a 
covering of $N^3$ of  degree 1 and the number 
of factors $S^2\times S^2$  is $(\sum_{i=1}^Dr_i)+D-1$. 
Thus any fiber sum decomposition of $M^4$ leads to a fiber sum whose associated 
covering is of degree 1. Then Lemma \ref{unique} implies that  the value of 
$(\sum_{i=1}^Dr_i)+D-1$ 
is independent on the fiber sum decomposition and thus the number of critical points, namely 
$2r+2D$, equals $\varphi(M,N)$,  as claimed.

\begin{remark}\label{iterate}
If $W^{2k}\to N^{k+1}$ is a spherical fibration 
(with connected fiber) then  it seems that the fiber sum of 
$W^{2k}$ with  $D$ spherical blocks $\#_{r_i}S^k\times S^k$ along 
$D$ fibers  is actually  diffeomorphic to the fiber sum of $W^{2k}$ 
with a single spherical block 
$\#_{(\sum_{i=1}^Dr_i)+D-1}S^k\times S^k$. 
This is not anymore clear if the fiber is not connected. 
\end{remark}

Consider now the case when $N^3$ is reducible. Thus either 
$\pi_1(N^3)=\Z$  or else $N^3$ splits as a connected sum 
$N_1^3\# N_2^3\#\cdots \# N_p^3$, where 
$N_i^3$ are either irreducible or  have fundamental group $\Z$. 
It is known that a closed orientable 3-manifold 
with fundamental group $\Z$ is diffeomorphic to $S^2\times S^1$, 
modulo the Poincar\'e Conjecture (see e.g. \cite{Hempel}).  
Let now $\widehat{N^3}$ be a finite covering of $N$  
whose fundamental   group is isomorphic $\pi_1(N)$. 
Then $\pi_1(\widehat{N^3})$ is either $\Z$ or splits 
as $\pi_1(N_1)*\pi_1(N^2)*\cdots *\pi_1(N_p)$. In the first case 
the manifolds are both $S^1\times S^2$ and thus they are diffeomorphic. 

Assume that the second alternative holds. 
By the affirmative solution of Kneser's Conjecture (see e.g. \cite{Wh,Ja}) we have 
that $\widehat{N^3}$ geometrically splits as $M_1^3\# M_2^3\#\cdots \# M_p^3$
where $M_j^3$ are closed 3-manifolds with $\pi_1(M_j^3)\cong \pi_1(N_j^3)$. 

\begin{lemma}
The irreducible 3-manifolds with the same 
fundamental group $M_j^3$ and $N_j^3$ are diffeomorphic. 
\end{lemma}
\begin{proof}
By the Jaco-Shalen-Johansson decomposition Theorem 
there exists a  family  (possibly empty) 
of incompressible tori in each one of the manifolds $M_j$ and $N_j$ such that 
each connected component of the complement is atoroidal. 

The family of tori is non-empty if there exists a $\Z\oplus \Z$ 
in the fundamental group. Therefore if the family of tori associated 
to one manifold is non-empty  the family of tori associated to the 
other one is also non-empty and both 3-manifolds are Haken, since they 
are irreducible and contain incompressible surfaces. 
Waldhausen's Theorem claims that a closed 3-manifold which 
is homotopy equivalent to a Haken manifold is actually diffeomorphic to it. 
Since both manifolds are irreducible they are aspherical and the isomorphism 
of their fundamental groups implies that they are homotopy equivalent. 

If the family of tori is empty the manifolds are atoroidal. 
In this case Thurston's geometrization Conjecture states that the manifolds 
are hyperbolic. Then, by Mostow's rigidity, two hyperbolic manifolds 
having the same fundamental group are isometric and hence diffeomorphic. 
\end{proof}
This implies that $\widehat{N^3}$ is diffeomorphic to $N^3$. 
Now we can apply the argument used in the case of non-co-Hopfian groups 
of irreducible 3-manifolds. This settles Proposition \ref{noncohopf}. 
\end{proof}

\bibliographystyle{amsplain}

\begin{thebibliography}{10}

\bibitem {AF1} D.Andrica and L.Funar, 
\textit{On smooth maps with finitely many critical
points}, J. London Math. Soc. 69(2004), 783--800,
{\it Addendum} 73(2006), 231--236.



\bibitem{AF2}
D.Andrica, L.Funar and E.Kudryavtseva, 
{\it On the minimal number of critical points of maps between 
closed manifolds},  Russian Journal of Mathematical Physics, special issue 
``Conference for the 65-th birthday of Nicolae Teleman'', 
(J.-P. Brasslet, A. Legrand, R. Longo, A. Mishchenko, Editors), 
16(2009), 363-370. 

\bibitem {Ant1} P.L.Antonelli, \textit{Structure theory for Montgomery-Samelson between
manifolds}, I, II Canad. J. Math. 21(1969), 170--179,  180--186.


\bibitem {Ant3} P.L.Antonelli, \textit{Differentiable Montgomery-Samelson fiberings with
finite singular sets}, Canad. J. Math. 21(1969),
1489--1495.

\bibitem{AT}
R. Ara\'ujo dos Santos, M.Tib\u{a}r, {\em Real map germs and higher open books}, arXiv:0801.3328. 


\bibitem{ArF}
E.Artin and R.Fox, 
\textit{Some wild cells and spheres in three-dimensional space}, 
Ann. of Math. 49(1948), 979-990. 



\bibitem{BBB}
L.Bessi\`eres, G.Besson, M.Boileau, S.Maillot and J.Porti, 
{\em Geometrisation of 3-manifolds}, 2009, to appear. 

\bibitem{Br}
Edward M.Brown, {\em Unknotted solid tori and genus one Whitehead manifolds},  
Trans. Amer. Math. Soc.  333(1992),  835--847. 

\bibitem{BLR}
D.Burghelea, R.Lashof and M.Rothenberg, {Groups of automorphisms of manifolds}, 
{\em Lecture Notes Math.} 473, Springer-Verlag, 1975. 

\bibitem{Cerf}
J.Cerf, Sur les diff\'eomorphismes de la sph\`ere de dimension trois 
$(\Gamma \sb{4}=0)$, 
{\em Lecture Notes in Math.}, 53, Springer, Berlin, 1968. 


\bibitem{CL}
P.T.Church and K.Lamotke, 
{\em Non-trivial polynomial isolated singularities}, 
Indag. Math. 37(1975), 149--154. 

\bibitem{CT1}
P.T.Church and J.G.Timourian,
{\em Differentiable maps with 0-dimensional critical set I},
Pacific J. Math. 41(1972), 615-630.

\bibitem{CT3}
 P.T.Church and J.G.Timourian,
{\em Continuous maps with 0-dimensional branch set},
Indiana Univ. Math. J. 23(1973/1974), 949-958.



\bibitem{CT2}
 P.T.Church and J.G.Timourian,
{\em Differentiable maps with 0-dimensional critical set II},
Indiana Univ. Math. J. 24(1974), 17-28.



\bibitem{Con}
P.E.Conner, {\it On the impossibility of fibring a certain manifold by a compact fibre}, Michigan Math.J. 5(1957), 249-255.


\bibitem{Dav}
R.Daverman,  Decompositions of manifolds, Academic Press, 1986. 

\bibitem {Dim} A.Dimca, Singularities and topology of
hypersurfaces, Springer-Verlag, Berlin, 1992.



\bibitem{DW}
A.Dold and H.Whitney, {\em Classification of oriented 
bundles over a 4-complex}, Ann. of Math. 69(1959), 667-677. 

\bibitem{EE}
C.Earle and J.Eells Jr.,  
{\em A fibre bundle description of Teichm\"uller theory}, 
  J. Differential Geometry  3(1969), 19--43. 



\bibitem{ES}
C.Earle and  J.H. Schatz, {\em Teichm\"uller theory for surfaces 
with boundary}, J. Differential Geometry 4(1970) , 169--185. 


\bibitem{FH}
E.R.Fadell and S.Y.Husseini, {Geometry and topology of configuration spaces}, 
{\em Springer Monograph Math.}, 2001.


\bibitem{FPZ}
L.Funar, C.Pintea and P.Zhang, {\em Examples of smooth maps with finitely
 many critical points in dimensions $(4,3)$, $(8,5)$ and $(16,9)$}, 
Proc. Amer. Math. Soc. 138(2010), 355-365. 


\bibitem {Hae1} A.Haefliger, \textit{Differentiable imbeddings}, Bull. Amer. Math. Soc. 67(1961), 109--112.

\bibitem{Hae2a}
A.Haefliger, {\em  
Plongements de vari\'et\'es dans le domaine stable}, 1964, S\'eminaire 
Bourbaki 1962/1963, 
Vol. 8,  Exp. No. 245, 63--77, Soc. Math. France, Paris, 1995. 

\bibitem{Hae2b}
A.Haefliger, {\em  
Plongements diff\'erentiables dans le domaine stable}, 
Comment. Math. Helv. 37(1962/1963), 155--176. 

\bibitem {Hae3} A.Haefliger, \textit{Differentiable embeddings of $S^{n}$ in $S^{n+q}$ for
$q>2$}, Ann. of Math. 83(1966), 402--436.


\bibitem{Ham}
M.-E.Hamstrom, 
{\em Regular mappings and the space of homeomorphisms on a $3$-manifold}, 
Mem. Amer. Math. Soc. No. 40(1961), 42p. 


\bibitem{H}
A.Hatcher, \textit{ Homeomorphisms of sufficiently large $P\sp{2}$-irreducible $3$-manifolds}, 
Topology  15(1976), 343--347.

\bibitem {Hat} 
A.Hatcher, \textit{A proof of a Smale conjecture, ${\rm Diff}(S\sp{3})\simeq {\rm O}(4)$}, 
Ann. of Math. (2) 117(1983), 553--607. 

\bibitem{HaM}
A.Hatcher and D.McCullough, {\em Finiteness of classifying spaces of relative diffeomorphism groups of $3$-manifolds},  Geom. Topol.  1(1997), 91--109. 





\bibitem{Hempel}
J.Hempel, 3-manifolds,  
reprint of the 1976 original,  
AMS Chelsea Publishing, 2004. 



\bibitem{HM}
W.Huebsch  and M.Morse, {\em Schoenflies extensions without interior
differential singularities}, Ann. of  Math.  76(1962), 18--54.


\bibitem{Ja}
W.Jaco, {\em Three manifolds with fundamental group a free product} 
Bull. Amer. Math. Soc. 75(1969), 972--977.



\bibitem{JW}
I.M.James and J.H.C.Whithead, {\it The homotopy theory of sphere bundles over spheres. II}, Proc. London Math. Soc. (3) 5(1955), 148--166. 



\bibitem{KN}
L.Kauffman and W.Neumann, 
{\em Products of knots, branched fibrations and sums of singularities}, 
Topology 16(1977),  369--393. 

\bibitem{King}
H.C.King,  {\em Topological type of isolated singularities}, 
Annals of Math. 107(1978), 385-397. 

\bibitem{King2}
H.C.King,  {\em Topology of isolated critical points of 
functions on singular spaces},  Stratifications, singularities and differential equations, II (Marseille, 1990; Honolulu, HI, 1990),  63--72, ((David Trotman and Leslie Charles Wilson, editors),  Travaux en Cours, 55, Hermann, Paris, 1997. 

\bibitem {Kos1} A.Kosinski, \textit{On the inertia group of
$\pi$- manifolds}, Amer. J. Math. 89(1967), 227--248.

\bibitem {Kos2} A.Kosinski, \textit{Differentiable manifolds}, Academic Press, London,
1993.


\bibitem{Lo}
E.Looijenga, {\em A note on polynomial isolated singularities}, 
Indag. Math.  33(1971),  418--421. 





\bibitem{Mil0}
J.Milnor, \textit{ On manifolds homeomorphic to the $7$-sphere},
  Ann. of Math. (2)  64(1956), 399--405. 

\bibitem {Mil1} J.Milnor, \textit{Differentiable structures on
spheres}, Amer. J. Math. 81(1959), 962--972.

\bibitem{Mil2}
J.Milnor, \textit{
A unique decomposition theorem for $3$-manifolds},  
Amer. J. Math.  84(1962), 1--7.

\bibitem{Milnor}
J.Milnor, Singular points of complex hypersurfaces, 
Princeton University Press, 1968. 

\bibitem {MonSam} 
D.Montgomery and H.Samelson, \textit{Fiberings with singularities},
Duke Math. J. 13(1946), 51--56.




\bibitem{PR}
P.E.Pushkar and Yu.B.Rudyak, {\em On the minimal number of critical points of functions on $h$-cobordisms},  Math. Res. Lett.  9(2002),  241--246.


\bibitem{Rees}
E.G.Rees, {\em On a question of Milnor concerning singularities of maps}, 
 Proc. Edinburgh Math. Soc. (2)  43(2000), 149--153. 

\bibitem{Ro}
D.Rolfsen, Knots and links, corrected reprint of the 1976 original, 
Mathematics Lecture Series, 7. Publish or Perish, Inc., 1990.




\bibitem{Rud}
L.Rudolph, {\em Isolated critical points of mappings from $R\sp 4$ to $R\sp 2$ and a natural splitting of the Milnor number of a classical fibered link. 
I. Basic theory; examples},  Comment. Math. Helv.  62(1987),  630--645. 

\bibitem{Sch}
R.Schultz, {\em
On the inertia group of a product of spheres},
Trans. Amer. Math. Soc. 156(1971), 137--153.

\bibitem{Sc}
P.Scott, {\em There are no fake Seifert fibre spaces with infinite $\pi \sb{1}$},   
Ann. of Math. (2)  117  (1983), no. 1, 35--70.

\bibitem{Sie}
L.C.Siebenmann, {\em Approximating cellular maps by homeomorphisms}, 
Topology 11(1972), 271-294. 

\bibitem{Sie2}
L.C.Siebenmann, {\em Deformation of homeomorphisms of stratified sets}, 
Comment. Math. Helv. 47(1972), 123-163. 

\bibitem{Sp}
E.H.Spanier, Algebraic topology, Springer-Verlag, New York-Berlin, 1981. 


\bibitem{St}
N.Steenrod, The topology of fibre bundles, 
reprint of the 1957 edition, Princeton University Press, 1999.


\bibitem{Takens}
F.Takens, {\em 
Isolated critical points of $C\sp{\infty }$ and $C\sp{\omega }$ functions}, 
Indag. Math. 29(1967), 238--243. 

\bibitem{Takens2}
F.Takens, {\em 
The minimal number of critical points of a function on a compact manifold and the Lusternik-Schnirelman category},  
Invent. Math.  6(1968), 197--244.


\bibitem {Tim} J.G.Timourian, \textit{Fiber bundles with discrete singular set},
J. Math. Mech. 18(1968), 61--70.



\bibitem{Wa}
F.Waldhausen, {\em On irreducible $3$-manifolds which are sufficiently large}, 
  Ann. of Math. (2)  87(1968), 56--88. 

\bibitem{WW}
Shicheng Wang and  Ying Qing Wu,  
{\em Covering invariants and co-Hopficity of $3$-manifold groups}, 
  Proc. London Math. Soc. (3)  68(1994),  203--224.

\bibitem{Wh}
J.H.C.Whitehead, {\em On finite cocycles and the sphere theorem}, 
Colloq. Math.  6(1958), 271--281. 

\bibitem{Woo}
L.M.Woodward, {\em The classification of orientable vector bundles 
over CW-complexes of small dimension}, Proc. Roy. Soc. Edinburgh, 
92A(1982), 175-179. 


\bibitem{YW}
Yu Fengchun and Shicheng Wang, {\em 
Covering degrees are determined by graph manifolds involved}, 
Comment. Math. Helv. 74(1999), 238--247. 


\end{thebibliography}

\end{document}